\documentclass{article}

\usepackage[a4paper]{geometry}
\usepackage[utf8]{inputenc}
\usepackage{amsfonts,amsmath,amssymb,amsthm,bm}
\usepackage{algorithm,algorithmic}
\usepackage{graphicx}
\usepackage[colorlinks=false, pdfborder={0 0 0}]{hyperref}
\usepackage{authblk}
\usepackage{setspace}
\usepackage{multirow}
\usepackage{xcolor}

\graphicspath{{fig/}}
\DeclareGraphicsExtensions{.pdf}

\theoremstyle{plain}

\newtheorem{theorem}{Theorem}[section]
\newtheorem{assumption}[theorem]{Assumption}
\newtheorem{lemma}[theorem]{Lemma}

\theoremstyle{remark}
\newtheorem{remark}{Remark}[section]

\numberwithin{equation}{section}

\DeclareMathOperator{\cond}{cond}

\DeclareMathOperator{\diag}{diag}

\DeclareMathOperator{\fl}{f{}l}

\DeclareMathOperator{\rank}{rank}

\newcommand*{\abs}[1]{\lvert#1\rvert}
\newcommand*{\bigabs}[1]{\big\lvert#1\big\rvert}
\newcommand*{\Bignorm}[1]{\Big\lVert#1\Big\rVert_{2}}
\newcommand*{\bignormfro}[1]{\big\lVert#1\big\rVert_{\fro}}
\newcommand*{\Bignormfro}[1]{\Big\lVert#1\Big\rVert_{\fro}}
\newcommand*{\bignorminf}[1]{\big\lVert#1\big\rVert_{\infty}}
\newcommand*{\epsgemv}{\epsilon_{\mathtt{gemv}}}
\newcommand*{\epsqr}{\epsilon_{\mathtt{qr}}}
\newcommand*{\epstrsv}{\epsilon_{\mathtt{trsv}}}
\newcommand*{\epsunmqr}{\epsilon_{\mathtt{unmqr}}}
\newcommand*{\fro}{\mathsf F}
\newcommand*{\herm}{^{*}}
\newcommand*{\iherm}{^{-*}}

\newcommand*{\machepsf}{\bm{u}_{f}}
\newcommand*{\machepso}{\bm{u}}
\newcommand*{\machepsr}{\bm{u}_{r}}
\newcommand*{\machepss}{\bm{u}_{s}}

\newcommand*{\norm}[1]{\lVert#1\rVert_{2}}
\newcommand*{\normbig}[1]{\left\lVert#1\right\rVert_{2}}
\newcommand*{\normfro}[1]{\lVert#1\rVert_{\fro}}
\newcommand*{\norminf}[1]{\lVert#1\rVert_{\infty}}
\newcommand*{\pinv}{^{\dagger}}

\newenvironment{keywords}{\medskip\textbf{Keywords:}}{}
\newenvironment{AMS}{\medskip\textbf{AMS subject classifications (2020).}}{}

\title{Mixed precision iterative refinement for least squares with
linear equality constraints and generalized least squares problems}
\author[1]{Bowen Gao}
\author[2]{Yuxin Ma}
\author[1,3]{Meiyue Shao}
\affil[1]{\normalsize School of Data Science,
Fudan University, Shanghai 200433, China

\texttt{Email: bwgao22@m.fudan.edu.cn, myshao@fudan.edu.cn}}
\affil[2]{Department of Numerical Mathematics,
Faculty of Mathematics and Physics, Charles University,
Sokolovsk\'{a} 49/83, 186 75 Praha 8, Czech Republic

\texttt{Email: yuxin.ma@matfyz.cuni.cz}}
\affil[3]{Shanghai Key Laboratory for Contemporary Applied Mathematics,
Fudan University, Shanghai 200433, China}

\begin{document}
\maketitle
\begin{abstract}
Recent development on mixed precision techniques has largely enhanced
the performance of various linear algebra solvers, one of which being
the solver for the least squares problem \(\min_{x}\norm{b-Ax}\).
By transforming least squares problems into augmented linear systems,
mixed precision techniques are capable of refining
the lower precision solution to the working precision.
In this paper, we propose mixed precision iterative refinement algorithms
for two variants of least squares problems---the least squares
problem with linear equality constraints~(LSE)
and the generalized least squares problem (GLS).
Both classical and GMRES-based iterative refinement can be applied to augmented
systems of these two problems to improve the accuracy of the solution.
For reasonably well-conditioned problems, our algorithms reduce
the execution time by a factor of \(40\%\) on average compared to
the fixed precision ones from LAPACK on the x86-64 architecture.
\end{abstract}


\section{Introduction}
\label{sec:introduction}
The least squares problem is one of the most
important topics in numerical linear algebra.
The least squares family contains many of the fundamental
problems encountered in physics~\cite{CZCY2023},
engineering~\cite{DH1972}, computer vision~\cite{HZ2004},
and other industrial applications~\cite{Bjorck1996,LH1995}.
Among those, the most important ones include the least squares with linear
equality constraints, the generalized least squares, the weighted least
squares, the rank-deficient least squares, and the total least squares, etc.
This work considers \emph{least squares problems with linear equality
constraints (LSE)} and \emph{generalized least squares problems (GLS)}.

\begin{enumerate}
\item  The LSE problem: Given matrices \(A\in\mathbb{C}^{m\times n}\)
and \(B\in\mathbb{C}^{p\times n}\) with \(p\leq n\leq m+p\)
and vectors \(b\in\mathbb{C}^{m}\) and \(d\in\mathbb{C}^{p}\),
the LSE problem is described as
\begin{equation}
\label{problem:lse}
\min_{x\in\mathbb{C}^{n}}\norm{Ax-b}\qquad\text{s.t.}\qquad Bx=d.
\end{equation}
Two requirements should be noted for the LSE problem.
\begin{enumerate}
\item  The linear constraint \(Bx=d\) must be \emph{consistent} for
any \(d\in\mathbb{C}^{p}\). This is equivalent to \(\rank(B)=p\).
\item  The solution \(x\) is \emph{unique} for any \(b\in\mathbb{C}^{m}\).
This is satisfied if and only if \(\mathcal{N}(A)\cap\mathcal{N}(B)=\{0\}\),
where \(\mathcal{N}(\cdot)\) represents the nullspace of a matrix,
i.e., \(\mathcal{N}(A)=\{x\colon Ax=0\}\).
This is equivalent to
\[
\rank\bigl(\begin{bmatrix}A \\ B\end{bmatrix}\bigr)=n.
\]
\end{enumerate}
In this paper, we shall only discuss LSE problems
that satisfy these two requirements, i.e.,
\begin{equation}
\label{assump:lse}
\rank(B)=p,\qquad\rank\bigl(\begin{bmatrix}A \\ B\end{bmatrix}\bigr)=n.
\end{equation}

The LSE problem is particularly important in applied regression analysis,
one of the universal methods of mathematical modelling.
Note that normal regression analysis takes the form \(b=Ax+e\),
which transforms into the least squares problem \(\min_{\beta}\norm{Ax-b}\).
If there exist certain equality constraints on \(x\),
e.g., from an exact observation or a physical law, then the standard
least squares problem becomes an LSE problem~\eqref{problem:lse}.
Therefore it is frequently encountered in applications such as
constrained surface fitting, constrained optimization, geodetic least
squares adjustment, signal processing, and other fields of interest.
See~\cite{Bjorck1996,GV2013} for more details.

\item  The GLS problem: Given matrices \(W\in\mathbb{C}^{n\times m}\)
and \(V\in\mathbb{C}^{n\times p}\) with \(m\leq n\leq m+p\) and
vector \(d\in\mathbb{C}^{n}\), the GLS problem is described as
\begin{equation}
\label{problem:gls}
\min_{x\in\mathbb{C}^{m},\,y\in\mathbb{C}^{p}}
\norm{y}\qquad\text{s.t.}\qquad Wx+Vy=d.
\end{equation}

Note that~\eqref{problem:gls} can be viewed as an LSE problem of the form
\begin{equation}
\label{eq:link}
\min_{x\in\mathbb{C}^m,\,y\in\mathbb{C}^p}\Bignorm{[0,I_{p}]
\begin{bmatrix}x \\ y\end{bmatrix}}\qquad\text{s.t.}
\qquad[W,V]\begin{bmatrix}x \\ y\end{bmatrix}=d.
\end{equation}
As discussed above, when
\[
\rank\bigl([W,V]\bigr)=n,\qquad
\rank\Bigl(\begin{bmatrix}0 & I_{p} \\ W & V\end{bmatrix}\Bigr)=m+p,
\]
or equivalently,
\begin{equation}
\label{assump:gls}
\rank\big([W,V]\big)=n,\qquad\rank(W)=m,
\end{equation}
the constraint is always consistent, and there exists a unique solution \(x\)
and a minimal 2-norm solution \(y\) to the system~\eqref{problem:gls}.

The GLS problem is crucial in economics and engineering as it is directly
associated with the generalized linear regression model (GLM).
GLM is a generalization of the standard linear regression model.
Given a matrix \(W\in\mathbb{C}^{n\times m}\) with full column rank,
the GLM is presented as
\begin{equation}
\label{eq:glm}
d=Wx+e,
\end{equation}
where \(e\) is a random vector with mean \(0\)
and covariance matrix \(\sigma^{2}Y\).
If \(\rank(Y)=p\), then \(Y\) can be factorized to \(Y=VV\herm\),
where \(V\in\mathbb{C}^{n\times p}\) has linearly independent columns.
Thus~\eqref{eq:glm} can be rewritten as
\[
d=Wx+Vy,
\]
where \(y\) is a random vector with mean \(0\)
and covariance matrix \(\sigma^{2}I\).
The GLS problem~\eqref{problem:gls} aims to minimize \(\sigma^{2}\).
We refer readers to~\cite{Paige1979-1,Paige1979-2} for more details.
\end{enumerate}

Notice from~\eqref{eq:link} that LSE and GLS problems have a similar structure.
As we shall show later in Section~\ref{sec:preliminary}, both problems can be
solved via the generalized RQ (GRQ) or generalized QR (GQR) factorization.
Therefore, these two problems are studied together in our work.

Mixed precision algorithms have attracted great attention
in recent years partly due to the evolution of machine
learning techniques and hardware infrastructure.
The industry is not satisfied with performing all floating-point
operations in a fixed precision (a.k.a.\ working precision),
as lower precision arithmetic conducts faster floating-point
operations and requires less data transfer.
However, exploiting lower precision arithmetic
may reduce the accuracy of the final solution.
A clever use of lower precision arithmetic offers the possibility
of greatly accelerating the computation, while still delivering
an accurate solution in the \emph{working precision} at the end.
This paves the way for mixed precision algorithms in various applications;
see the surveys~\cite{Survey2021,HM2022}.

Existing works on mixed precision algorithms in numerical
linear algebra have already covered a number of topics.
Those include algorithms for solving linear
systems~\cite{CH2017,CH2018,LLLKD2006} as well as algorithms
for solving least squares problems~\cite{CHP2020}.
A brief investigation into half precision (\texttt{FP16})
arithmetic is discussed in~\cite{HTDH2018,HWTD2017}.
Topics related to eigenvalues and singular values are more complicated.
Still, there exist several studies on this matter,
including the symmetric eigensolvers~\cite{GMS2025,KMS2023,OA2018,OA2019,
OA2020,XLLSR2024} and nonsymmetric eigensolvers~\cite{BKS2023}.

Let us focus more closely on linear systems and least squares problems.
The modern mixed precision approach for solving
linear systems has been adequately presented by Carson and Higham
in~\cite{CH2018}, which calculates the most expensive part---LU
factorization---in lower precision and uses the classical iterative
refinement algorithm to refine the solution to the working precision.
For ill-conditioned matrices, however, this approach may fail to converge.
Therefore~\cite{CH2018} applies a GMRES solver,
capable of utilizing the lower precision solution as the initial
guess and lower precision LU factors as preconditioner,
in such scenarios to further expand the feasibility of the algorithm.
Since least squares problems can be transformed to 
certain linear systems, mixed precision techniques
for linear systems carry over to least squares problems.
Carson, Higham, and Pranesh proposed a mixed precision GMRES-based
iterative refinement algorithm by performing a lower precision
QR factorization to compute the initial guess and then reusing
the QR factors as preconditioner to solve least squares
problems through augmented systems~\cite{CHP2020,HP2021}.
Mixed precision algorithms on weighted least
squares~\cite{CO2025} and total least squares
problems~\cite{OC2024} have been studied by Carson and Oktay.

In this paper, we propose mixed precision classical
and GMRES-based iterative refinement algorithms for LSE
and GLS problems by utilizing the GRQ or GQR factorization.
We analyse the four-precision iterative refinement algorithm
and prove that the mixed precision algorithm is able to refine
the lower precision solutions to the working precision level.
We also show that when employing a higher precision in
the computation of residuals, the iterative refinement
method is able to improve the accuracy of the solution.
Our numerical experiments further indicate that for well-conditioned problems,
these mixed precision algorithms offer significant speedup compared
to the fixed precision algorithms without sacrificing accuracy.

The rest of the paper is structured as follows.
In Section~\ref{sec:preliminary} we describe standard (fixed precision)
approaches for solving LSE and GLS problems using the GRQ and GQR factorization.
Based on the augmented systems, we propose a mixed precision iterative
refinement algorithm for each problem in Section~\ref{sec:algorithm}.
We then present our four-precision rounding
error analysis in Section~\ref{sec:accuracy}.
In Section~\ref{sec:gmres} we identify the limitations of the classical
iterative refinement algorithm and design preconditioners for the GMRES
solver to overcome this issue for ill-conditioned matrices.
Numerical experiments are provided in Section~\ref{sec:experiments} to
demonstrate the accuracy and efficiency of the mixed precision algorithms.

We first introduce some notation used throughout the paper.
We employ the MATLAB colon notation to present the submatrices and subvectors.
Note that \(\norm{\cdot}\) indicates the \(2\)-norm,
\(\norminf{\cdot}\) indicates the infinity norm,
and \(\normfro{\cdot}\) is utilized for the Frobenius norm.
Let \(\kappa_{2}(A)=\norm{A}\cdot\norm{A\pinv}\) and 
\(\kappa_{\infty}(A)=\norminf{A}\cdot\norminf{A\pinv}\)
denote the matrix condition numbers under different norms,
where \(A\pinv\) represents the Moore--Penrose pseudoinverse of \(A\).
Additionally, \(\abs{\cdot}\) denotes the absolute value of
each matrix element, and \(\abs{A_{i,j}}\leq\abs{B_{i,j}}\)
for all \(i\) and \(j\) is simplified as \(\abs{A}\leq\abs{B}\).
The notation \(x\lesssim y\) signifies that \(x\) is
less than or equal to \(y\) up to a multiplicative constant,
i.e., \(x\leq C\cdot y\) for some constant \(C>0\).
\section{Algorithms for LSE and GLS problems}
\label{sec:preliminary}
In this section, we review the existing algorithms based on the GRQ or GQR
factorization and corresponding augmented systems for LSE and GLS problems.
We remark that the notation in LSE and GLS are \emph{independent},
even if some variables may have the same name.

\subsection{Algorithms for the LSE problem}
\label{subsec:preliminary-lse}
For the LSE problem~\eqref{problem:lse},
we first present the nullspace approach via the GRQ factorization
proposed by~\cite{ABD1992, Bjorck1996} as follows:

\begin{enumerate}
\item  Compute the GRQ factorization of \(B\) and \(A\) such that
\[
B=[0,R]Q,\qquad A=ZTQ,
\]
where \(R\in\mathbb{C}^{p\times p}\) is an upper triangular matrix, \(T\in\mathbb{C}^{m\times n}\)
is an upper rectangular matrix, and \(Z\in\mathbb{C}^{m\times m}\),
\(Q\in\mathbb{C}^{n\times n}\) are unitary matrices.
\item  Compute \(y_{2}\in\mathbb{C}^{p}\) by
solving the triangular system \(Ry_{2}=d\).
\item  Compute \(y_{1}\in\mathbb{C}^{n-p}\) by solving the triangular
system \(T_{11}y_{1}=c_{1}-T_{12}y_{2}\) with the partition of
\[
T=\begin{bmatrix}T_{11} & T_{12} \\ 0 & T_{22}\end{bmatrix}\qquad
\text{and}\qquad c=Z\herm b=\begin{bmatrix}c_{1} \\ c_{2}\end{bmatrix}.
\]
\item  Assign \(x\gets Q\herm\begin{bmatrix}y_{1} \\ y_{2}\end{bmatrix}\).
\end{enumerate}

The augmented (or saddle-point) system leads to another
approach to solving the LSE problem~\cite{ST2022}.
For instance, a straightforward approach solves
the LSE problem via the augmented system of
\[
\begin{bmatrix}A\herm A & B\herm \\ B & 0\end{bmatrix}
\begin{bmatrix}x \\ v\end{bmatrix}=\begin{bmatrix}A\herm b \\ d\end{bmatrix}.
\]
A 3-block version is expressed as
\begin{equation}
\label{eq:aug-LSE}
\begin{bmatrix}I_{m} & 0 & A \\ 0 & 0 & B \\ A\herm & B\herm & 0
\end{bmatrix}\begin{bmatrix}r \\ -v \\ x\end{bmatrix}
=\begin{bmatrix}b \\ d \\ 0\end{bmatrix},
\end{equation}
which circumvents forming \(A\herm A\) explicitly.
In Section~\ref{sec:algorithm}, the augmented system~\eqref{eq:aug-LSE} will
be used to perform the iterative refinement process for the LSE problem.
For further discussion on saddle-point systems,
see~\cite{BGL2005,Rozloznik2018}.

\subsection{Algorithms for the GLS problem}
\label{subsec:preliminary-gls}
As shown in Section~\ref{subsec:preliminary-lse}, the nullspace
approach utilizes the GRQ factorization to solve the LSE problem.
Similarly, Paige~\cite{Paige1979-1,Paige1979-2} proposed an algorithm
to address the GLS problem through the GQR factorization as follows:

\begin{enumerate}
\item  Compute the GQR factorization of \(W\) and \(V\) such that
\[
W=Q\begin{bmatrix}R \\ 0\end{bmatrix},\qquad V=QTZ\qquad\text{with}
\qquad T=\begin{bmatrix}T_{11} & T_{12} \\ 0 & T_{22}\end{bmatrix},
\]
where \(R\in\mathbb{C}^{m\times m}\), \(T_{22}\in\mathbb{C}^{(n-m)
\times(n-m)}\) are upper triangular matrices, and \(Q\in\mathbb{C}^
{n\times n}\), \(Z\in\mathbb{C}^{p\times p}\) are unitary matrices.
\item  Compute \(s_{2}\in\mathbb{C}^{n-m}\) by solving
the triangular system \(T_{22}s_{2}=c_{2}\),
where \(c=Q\herm d=\begin{bmatrix}c_{1} \\ c_{2}\end{bmatrix}\).
\item  Compute \(x\in\mathbb{C}^{m}\) by solving
the triangular system \(Rx=c_{1}-T_{12}s_{2}\).
\item  Let \(s_{1}=0\in\mathbb{C}^{m-n+p}\) and assign
\(y\gets Z\herm\begin{bmatrix}s_{1} \\ s_{2}\end{bmatrix}\).
\end{enumerate}

Similarly, the augmented system of the GLS problem is expressed as
\[
\begin{bmatrix}VV\herm & W \\ W\herm & 0\end{bmatrix}
\begin{bmatrix}z \\ x\end{bmatrix}=\begin{bmatrix}d \\ 0\end{bmatrix}.
\]
The GLS problem can also be solved via the 3-block version of
\begin{equation}
\label{eq:aug-GLS}
\begin{bmatrix} I_{p} & V\herm & 0 \\ V & 0 & W \\ 0 & W\herm & 0
\end{bmatrix}\begin{bmatrix}y \\ -z \\ x\end{bmatrix}
=\begin{bmatrix}0 \\ d \\ 0\end{bmatrix},
\end{equation}
which avoids forming \(VV\herm\) explicitly,
and can be used for iterative refinement.
\section{Mixed precision algorithms}
\label{sec:algorithm}
In this section, we propose mixed precision algorithms based
on iterative refinement to solve LSE and GLS problems.
The basic idea is to perform the more costly GRQ or GQR factorization
in lower precision and then refine the lower precision solution
by iterative refinement to reach the working precision level.

Let \(\machepso\) be the unit roundoff at the working precision.
Similarly to~\cite{CHP2020}, we employ four precisions \(1\gg\,
\machepsf\geq\machepss\geq\machepso\geq\machepsr\) in our analysis,
where \(\machepsf\), \(\machepsr\), and \(\machepss\) stand for
the unit roundoff used in the matrix factorization, the computation
of residuals, and the correction equation solver, respectively.
The general framework of mixed precision LSE/GLS algorithms is as follows:

\begin{enumerate}
\item  Compute the initial guess of the least squares problem
via matrix factorization at precision~\(\machepsf\).
The matrix factors are kept in the memory.
\item  Let \(\tilde{F}\) and \(\tilde{s}\) be
the augmented matrix and vector, respectively.
Perform iterative refinement at multiple precisions iteratively
to solve the augmented system \(\tilde{F}\tilde{u}=\tilde{s}\)
in order to reach the required accuracy as follows:
\begin{enumerate}
\item  Compute the residual \(\tilde{r}=\tilde{s}-\tilde{F}\tilde{u}\)
of the augmented system at precision \(\machepsr\).
\item  Solve the linear system \(\tilde{F}\Delta\tilde{u}=\tilde{r}\)
to obtain \(\Delta\tilde{u}\) at precision \(\machepss\).
By making use of the existing matrix factorization,
it suffices to solve several triangular linear systems.
\item  Update \(\tilde{u}\gets\tilde{u}+\Delta\tilde{u}\)
at precision \(\machepso\).
\end{enumerate}
\end{enumerate}

For LSE and GLS problems, as discussed in Section~\ref{sec:preliminary},
we shall use the nullspace approach and Paige's algorithm, respectively,
to perform the first step of the process at precision \(\machepsf\)
so as to obtain the initial guess and the GRQ/GQR factorization.
In this section, we further develop the four-precision
classical iterative refinement method.

\subsection{Classical iterative refinement for the LSE problem}
\label{subsec:mpalgo-lse}
We first focus on the LSE problem~\eqref{problem:lse}.
After obtaining the initial approximation \(x_{0}\) and GRQ
factorization of \((B,A)\) at precision \(\machepsf\), i.e.,
\begin{equation}
\label{eq:grq}
B=R_{B}Q,\qquad A=ZTQ,
\end{equation}
where \(R_{B}=[0,R]\); \(Q\in\mathbb{C}^{n\times n}\)
and \(Z\in\mathbb{C}^{m\times m}\) are unitary matrices,
\(R\in\mathbb{C}^{p\times p}\) is an upper triangular matrix,
\(T\in\mathbb{C}^{m\times n}\) is an upper rectangular matrix,
we compute the residual \(r_{0}=b-Ax_{0}\)
and \(v_{0}\) by solving the triangular system \(R\herm v_{0}
=(QA\herm r_{0})(n-p+1:n)\) at precision \(\machepsf\), respectively.

Given \(x\gets x_{0}\), \(r\gets r_{0}\), and \(v\gets v_{0}\),
each step of the classical iterative refinement
on~\eqref{eq:aug-LSE} is conducted as follows:

\begin{enumerate}
\item  Compute the residuals \(f_{1}\in\mathbb{C}^{m}\),
\(f_{2}\in\mathbb{C}^{p}\), and \(f_{3}\in\mathbb{C}^{n}\) of
the 3-block saddle-point system at precision \(\machepsr\):
\[
\begin{bmatrix}f_{1} \\ f_{2} \\ f_{3}\end{bmatrix}
=\begin{bmatrix}b \\ d \\ 0\end{bmatrix}-\begin{bmatrix}
I_{m} & 0 & A \\ 0 & 0 & B \\ A\herm & B\herm & 0\end{bmatrix}
\begin{bmatrix}r \\ -v \\ x\end{bmatrix}.
\]
\item  Compute \(\Delta r\in\mathbb{C}^{m}\), \(\Delta v\in\mathbb{C}^{p}\),
and \(\Delta x\in\mathbb{C}^{n}\) by solving the saddle-point
linear system at precision \(\machepss\):
\begin{equation}
\label{eq:correction-system-lse}
\begin{bmatrix} I_{m} & 0 & A \\ 0 & 0 & B \\ A\herm & B\herm & 0
\end{bmatrix}\begin{bmatrix}\Delta r \\ -\Delta v \\ \Delta x\end{bmatrix}
=\begin{bmatrix}f_{1} \\ f_{2} \\ f_{3}\end{bmatrix}.
\end{equation}
\item  Update \(r\), \(v\), and \(x\) at precision \(\machepso\):
\[
\begin{bmatrix}r \\ v \\ x\end{bmatrix}\gets\begin{bmatrix}r \\ v \\ x
\end{bmatrix}+\begin{bmatrix}\Delta r \\ \Delta v \\ \Delta x\end{bmatrix}.
\]
\end{enumerate}
This process is repeated iteratively until
the following stopping criteria are satisfied:
\begin{equation}
\label{eq:lse-stop}
\begin{split}
\norm{f_{1}} & \leq\mathtt{tol}\bigl(
\norm{b}+\norm{r}+\normfro{A}\norm{x}\bigr),\\
\norm{f_{2}} & \leq\mathtt{tol}\bigl(\norm{d}+\normfro{B}\norm{x}\bigr),\\
\norm{f_{3}} & \leq\mathtt{tol}\bigl(
\normfro{A}\norm{r}+\normfro{B}\norm{v}\bigr).
\end{split}
\end{equation}
The threshold \texttt{tol} is here chosen to be \(\mathcal{O}(\machepso)\).

In Step 2 of the procedure, we solve the linear
system~\eqref{eq:correction-system-lse} using the GRQ factors.
Substituting \(A\) and \(B\) into the GRQ factorization~\eqref{eq:grq},
the linear system~\eqref{eq:correction-system-lse} is rewritten as
\[
\begin{bmatrix}Z & 0 & 0 \\ 0 & I_{p} & 0 \\
0 & 0 & Q\herm\end{bmatrix}\begin{bmatrix} I_{m} & 0 & T \\
0 & 0 & R_{B} \\T\herm & R_{B}\herm & 0\end{bmatrix}
\begin{bmatrix}Z\herm & 0 & 0 \\0 & I_{p} & 0 \\ 0 & 0 & Q\end{bmatrix}
\begin{bmatrix}\Delta r \\ -\Delta v \\ \Delta x\end{bmatrix}
=\begin{bmatrix}f_{1} \\ f_{2} \\ f_{3}\end{bmatrix}
\]
with unknown variables \(\Delta r\), \(\Delta v\), and \(\Delta x\).
Let \(q=Z\herm\Delta r\in\mathbb{C}^{m}\) and \(y=Q\Delta x\in\mathbb{C}^{n}\).
Since \(Z\) and \(Q\) are unitary matrices, we have
\[
\begin{bmatrix} I_{m} & 0 & T \\ 0 & 0 & R_{B} \\
T\herm & R_{B}\herm & 0\end{bmatrix}\begin{bmatrix}q \\ -\Delta v \\ y
\end{bmatrix}=\begin{bmatrix}Z\herm f_{1} \\ f_{2} \\ Qf_{3}\end{bmatrix},
\]
which is equivalent to
\begin{equation}
\label{eq:secondary-lse}
\begin{split}
& q=Z\herm f_{1}-Ty,\qquad R_{B}y=f_{2},\qquad
T\herm q-\begin{bmatrix}0 \\ R\herm \Delta v\end{bmatrix}=Qf_{3}.
\end{split}
\end{equation}
Note that the unknown variables in~\eqref{eq:secondary-lse}
are \(q\), \(y\), and \(\Delta v\).
We may obtain \(\Delta r\) and \(\Delta x\) via
\(\Delta r=Zq\) and \(\Delta x=Q\herm y\), respectively.
Therefore, our goal shifts to solve~\eqref{eq:secondary-lse}.

Denote \(u=Qf_{3}\) and \(w=Z\herm f_{1}\), and let
\(q\), \(y\), \(u\), \(w\), and \(T\) be partitioned as
\[
q=\begin{bmatrix}q_{1} \\ q_{2} \end{bmatrix},\qquad
y=\begin{bmatrix}y_{1} \\ y_{2} \end{bmatrix},\qquad
u=\begin{bmatrix}u_{1} \\ u_{2}\end{bmatrix},\qquad
w=\begin{bmatrix}w_{1} \\ w_{2}\end{bmatrix},\qquad
T=\begin{bmatrix}T_{11} & T_{12} \\ 0 & T_{22}\end{bmatrix},
\]
with \(q_{1}\in\mathbb{C}^{n-p}\), \(q_{2}\in\mathbb{C}^{m-n+p}\),
\(y_{1}\in\mathbb{C}^{n-p}\), \(y_{2}\in\mathbb{C}^{p}\),
\(u_{1}\in\mathbb{C}^{n-p}\), \(u_{2}\in\mathbb{C}^{p}\),
\(w_{1}\in\mathbb{C}^{n-p}\), \(w_{2}\in\mathbb{C}^{m-n+p}\),
\(T_{11}\in\mathbb{C}^{(n-p)\times(n-p)}\), \(T_{12}\in\mathbb{C}^{(n-p)
\times p}\), and \(T_{22}\in\mathbb{C}^{(m-n+p)\times p}\).
Note that \(T_{11}\) and~\(R\) are both nonsingular
if the assumption~\eqref{assump:lse} is satisfied.
Then~\eqref{eq:secondary-lse} is reformulated as
\begin{subequations}
\label{eq:sum-lse}
\begin{align}
\label{eq:y1} & T_{11}y_{1}=w_{1}-q_{1}-T_{12}y_{2},\\
\label{eq:q2} & q_{2}=w_{2}-T_{22}y_{2},\\
\label{eq:y2} & Ry_{2}=f_{2},\\
\label{eq:q1} & T_{11}\herm q_{1}=u_{1},\\
\label{eq:deltav} & R\herm\Delta v=T_{12}\herm q_{1}+T_{22}\herm q_{2}-u_{2}.
\end{align}
\end{subequations}

To compute \(q\), \(y\), and \(\Delta v\), we start by solving
the triangular systems~\eqref{eq:y2} and~\eqref{eq:q1}
to obtain \(y_{2}\) and \(q_{1}\), respectively.
We then obtain \(y_{1}\) and \(q_{2}\) via~\eqref{eq:y1} and~\eqref{eq:q2}.
Finally, \(\Delta v\) is computed by solving
the triangular system~\eqref{eq:deltav}.
This procedure for solving the correction system
is outlined in Algorithm~\ref{alg:lse-solver}.
Combined with the mixed precision framework at the beginning of the section,
we summarize this mixed precision LSE algorithm in Algorithm~\ref{alg:mplse}.

\begin{algorithm}[!tb]
\caption{A correction system solver algorithm for the LSE problem}
\label{alg:lse-solver}
\begin{algorithmic}[1]
\REQUIRE  A GRQ factorization of \((B, A)\) with the upper triangular matrix
\(R\in\mathbb{C}^{p\times p}\), the upper rectangular \(T\in\mathbb{C}^{m\times n}\),
and unitary matrices \(Z\in\mathbb{C}^{m\times m}\),
\(Q\in\mathbb{C}^{n\times n}\); vectors \(f_{1}\in\mathbb{C}^{m}\),
\(f_{2}\in\mathbb{C}^{p}\), and \(f_{3}\in\mathbb{C}^{n}\).
\ENSURE  Vectors \(x_{1}\in\mathbb{C}^{m}\), \(x_{2}\in\mathbb{C}^{p}\),
and \(x_{3}\in\mathbb{C}^{n}\) satisfying
\[
\begin{bmatrix}I_{m} & 0 & A \\ 0 & 0 & B \\ A\herm & B\herm & 0
\end{bmatrix}\begin{bmatrix}x_{1} \\ -x_{2} \\ x_{3}\end{bmatrix}
=\begin{bmatrix}f_{1} \\ f_{2} \\ f_{3}\end{bmatrix}.
\]
\hspace{-2\algorithmicindent}
\textbf{Interface}: \((x_{1},x_{2},x_{3})=\mathtt{CorrectionSystemSolverLSE}
(R,T,Z,Q,f_{1},f_{2},f_{3})\)

\medskip
\STATE  \(u\gets Qf_{3}\), i.e., \(u_{1}=(Qf_{3})(1:n-p)\)
and \(u_{2}=(Qf_{3})(n-p+1:n)\).
\STATE  \(w\gets Z\herm f_{1}\), i.e., \(w_{1}=(Z\herm f_{1})(1:n-p)\)
and \(w_{2}=(Z\herm f_{1})(n-p+1:m)\).
\STATE  Compute \(y_{2}\) by solving the triangular system \(Ry_{2}=f_{2}\).
\STATE  Compute \(q_{1}\) by solving the triangular
system \(T_{11}\herm q_{1}=u_{1}\).
\STATE  Compute \(y_{1}\) by solving the triangular
system \(T_{11}y_{1}=w_{1}-q_{1}-T_{12}y_{2}\).
\STATE  \(q_{2}\gets w_{2}-T_{22}y_{2}\).
\STATE  \(x_{1}\gets Zq\) and \(x_{3}\gets Q\herm y\).
\STATE  Compute \(x_{2}\) by solving the triangular system
\(R\herm x_{2}=T_{12}\herm q_{1}+T_{22}\herm q_{2}-u_{2}\).
\RETURN  \((x_{1},x_{2},x_{3})\).
\end{algorithmic}
\end{algorithm}

\begin{algorithm}[!tb]
\caption{Mixed precision LSE algorithm}
\label{alg:mplse}
\begin{algorithmic}[1]
\REQUIRE  Matrices \(A\in\mathbb{C}^{m\times n}\),
\(B\in\mathbb{C}^{p\times n}\) and vectors \(b\in\mathbb{C}^{m}\),
\(d\in\mathbb{C}^{p}\) storing at precision \(\machepso\);
the maximal number of iterations \texttt{maxit} with a default value \(40\);
the tolerance parameter \texttt{tol} for the stopping criteria.
\ENSURE  A vector \(x\in\mathbb{C}^{n}\) stored at precision \(\machepso\)
approximating the solution of the LSE problem~\eqref{problem:lse}.

\medskip
\STATE  Compute the GRQ factorization of \((B,A)\) at
precision \(\machepsf\), i.e., \(B=[0,R]Q\), \(A=ZTQ\).
\STATE  Compute the initial guess \(x_{0}\) via the GRQ factorization at
precision \(\machepsf\), and store \(x\gets x_{0}\) at precision \(\machepso\).
\STATE  Compute \(r_{0}\gets b-Ax_{0}\) at precision \(\machepsf\),
and store \(r\gets r_{0}\) at precision \(\machepso\).
\STATE  Solve \(v_{0}\) by the triangular system
\(R\herm v_{0}=(Q\cdot A\herm r)(n-p+1:n)\) at precision
\(\machepsf\) with \(A\herm r\) computed at precision \(\machepsr\),
and store \(v\gets v_{0}\) at precision~\(\machepso\).
\FOR{\(i=1,\dotsc,\mathtt{maxit}\)}
\STATE  Compute the residuals \(f_{1}\gets b-r-Ax\), \(f_{2}\gets d-Bx\),
and \(f_{3}\gets-A\herm r+B\herm v\) at precision~\(\machepsr\),
and store them at precision \(\machepso\).
\IF{the stopping criteria~\eqref{eq:lse-stop} are satisfied}
\RETURN  \(x\).
\ENDIF
\STATE  Compute \(\Delta r\), \(\Delta v\), and \(\Delta x\)
at precision \(\machepss\) by function
\[
(\Delta r,\Delta v,\Delta x)=\mathtt{CorrectionSystemSolverLSE}
(R,T,Z,Q,f_{1},f_{2},f_{3})
\]
in Algorithm~\ref{alg:lse-solver}, and store them at precision \(\machepso\).
\STATE  Update and store \(r\gets r+\Delta r\), \(x\gets x+\Delta x\),
and \(v\gets v+\Delta v\) at precision \(\machepso\).
\ENDFOR
\end{algorithmic}
\end{algorithm}

Finally, we measure the complexity of the algorithm
in terms of the number of floating-point operations.
Let \(t_{f}\), \(t_{w}\), \(t_{r}\), and \(t_{s}\) represent the total number
of floating-point operations in Algorithm~\ref{alg:mplse} at precision
\(\machepsf\), \(\machepso\), \(\machepsr\), and \(\machepss\), respectively.
Then
\[
t_{f}=\mathcal{O}(mn^{2}+np^{2}).
\]
And in each refinement step,
\begin{align*}
t_{w} & =\mathcal{O}(m+n+p),\quad t_{r}=4(mn+np)+\mathcal{O}(m+n+p),\\
t_{s} & =2(mn+n^{2})+4mp+4(p^{2}+(n-p)^{2})+\mathcal{O}(m+n+p).
\end{align*}
In a common scenario where \(m\geq n\geq p\), we obtain
\[
t_{f}=\mathcal{O}(mn^{2}),\qquad t_{w}=\mathcal{O}(m),\qquad
t_{r}=\mathcal{O}(mn),\qquad t_{s}=\mathcal{O}(mn).
\]
It is apparent that \(t_{f}\) is the leading term
when the matrix dimension \(n\) is relatively large.

\subsection{Classical iterative refinement for the GLS problem}
\label{subsec:mpalgo-gls}
We derive a similar method for the GLS problem~\eqref{problem:gls}.
As discussed in Section~\ref{sec:preliminary},
we calculate the initial approximation \(x_{0}\), \(y_{0}\), and
the GQR factorization of \((W,V)\) at precision \(\machepsf\), i.e.,
\[
W=QR_{W},\qquad V=QTZ\qquad\text{with}\qquad
T=\begin{bmatrix}T_{11} & T_{12} \\ 0 & T_{22}\end{bmatrix}\quad
\text{and}\quad R_{W}=\begin{bmatrix}R \\ 0\end{bmatrix},
\]
where \(Q\in\mathbb{C}^{n\times n}\) and \(Z\in\mathbb{C}^{p\times p}\)
are unitary matrices, \(R\in\mathbb{C}^{m\times m}\) and
\(T_{22}\in\mathbb{C}^{(n-m)\times(n-m)}\) are upper triangular matrices.
We then solve the triangular system \(T_{22}\herm v=(Zy)(p-n+m+1:p)\)
and compute \(z_{0}\gets Q\herm\bar{v}\) at precision \(\machepsf\),
where \(\bar{v}=[0,v\herm]\herm\).

Given \(x\gets x_{0}\), \(y\gets y_{0}\), and \(z\gets z_{0}\),
the classical iterative refinement on the 3-block saddle-point
system~\eqref{eq:aug-GLS} is conducted as follows:

\begin{enumerate}
\item  Compute the residuals \(f_{1}\in\mathbb{C}^{p}\),
\(f_{2}\in\mathbb{C}^{n}\), and \(f_{3}\in\mathbb{C}^{m}\)
of the 3-block saddle-point system at precision \(\machepsr\):
\[
\begin{bmatrix}f_{1} \\ f_{2} \\ f_{3}\end{bmatrix}
=\begin{bmatrix}0 \\ d \\ 0\end{bmatrix}
-\begin{bmatrix} I_{p} & V\herm & 0 \\ V & 0 & W \\ 0 & W\herm & 0
\end{bmatrix}\begin{bmatrix}y \\ -z \\ x\end{bmatrix}.
\]
\item  Compute \(\Delta y\in\mathbb{C}^{p}\), \(\Delta z\in\mathbb{C}^{n}\),
and \(\Delta x\in\mathbb{C}^{m}\) by solving
the saddle-point linear system at precision \(\machepss\):
\begin{equation}
\label{eq:correction-system-gls}
\begin{bmatrix} I_{p} & V\herm & 0 \\ V & 0 & W \\ 0 & W\herm & 0
\end{bmatrix}\begin{bmatrix}\Delta y \\ -\Delta z \\ \Delta x\end{bmatrix}
=\begin{bmatrix}f_{1} \\ f_{2} \\ f_{3}\end{bmatrix}.
\end{equation}
\item  Update \(y\), \(z\), and \(x\) at precision \(\machepso\):
\[
\begin{bmatrix}y \\ z \\ x\end{bmatrix}\gets\begin{bmatrix}y \\ z \\ x
\end{bmatrix}+\begin{bmatrix}\Delta y \\ \Delta z \\ \Delta x\end{bmatrix}.
\]
\end{enumerate}
This process for the GLS problem concludes when the following stopping
criteria are satisfied, with \(\mathtt{tol} = \mathcal{O}(\machepso)\):
\begin{equation}
\label{eq:gls-stop}
\begin{split}
\norm{f_{1}} & \leq\mathtt{tol}\big(\norm{y}+\normfro{V}\norm{z}\big),\\
\norm{f_{2}} & \leq\mathtt{tol}\big(\norm{d}
+\normfro{W}\norm{x}+\normfro{V}\norm{y}\big),\\
\norm{f_{3}} & \leq\mathtt{tol}\normfro{W}\norm{z}.
\end{split}
\end{equation}

Denote \(u=Q\herm f_{2}\) and \(w=Zf_{1}\), and let \(g=Z\Delta y\),
\(h=Q\herm\Delta z\), \(u\), \(w\), and \(T\) be partitioned as
\[
g=\begin{bmatrix}g_{1} \\ g_{2}\end{bmatrix},\qquad
h=\begin{bmatrix}h_{1} \\ h_{2}\end{bmatrix},\qquad
u=\begin{bmatrix}u_{1} \\ u_{2}\end{bmatrix},\qquad
w=\begin{bmatrix}w_{1} \\ w_{2}\end{bmatrix},\qquad
T=\begin{bmatrix}T_{11} & T_{12} \\ 0 & T_{22}\end{bmatrix},
\]
with \(g_{1}\in\mathbb{C}^{p-n+m}\), \(g_{2}\in\mathbb{C}^{n-m}\),
\(h_{1}\in\mathbb{C}^{m}\), \(h_{2}\in\mathbb{C}^{n-m}\),
\(u_{1}\in\mathbb{C}^{m}\), \(u_{2}\in\mathbb{C}^{n-m}\),
\(w_{1}\in\mathbb{C}^{p-n+m}\), \(w_{2}\in\mathbb{C}^{n-m}\),
\(T_{11}\in\mathbb{C}^{m\times(p-n+m)}\), \(T_{12}\in\mathbb{C}
^{m\times(n-m)}\), and \(T_{22}\in\mathbb{C}^{(n-m)\times(n-m)}\).
Note that \(T_{22}\) and~\(R\) are both nonsingular
if the assumption~\eqref{assump:gls} is satisfied.
Then similarly to the procedure in Section~\ref{subsec:mpalgo-lse},
\eqref{eq:correction-system-gls} is reformulated as
\begin{equation}
\label{eq:gls-solver}
\begin{split}
& g_{1}=w_{1}-T_{11}\herm h_{1},\\
& T_{22}\herm h_{2}=w_{2}-g_{2}-T_{12}\herm h_{1},\\
& R\Delta x=u_{1}-T_{11}g_{1}-T_{12}g_{2},\\
& T_{22}g_{2}=u_{2},\\
& R\herm h_{1}=f_{3}.
\end{split}
\end{equation}
Note that~\eqref{eq:gls-solver} can be solved through
the same order as in Section~\ref{subsec:mpalgo-lse}.
This procedure is formally outlined in Algorithm~\ref{alg:gls-solver}.
Incorporating Algorithm~\ref{alg:gls-solver} with the mixed precision
framework discussed at the beginning of the section, we obtain
the mixed precision GLS algorithm as in Algorithm~\ref{alg:mpgls}.

\begin{algorithm}[!tb]
\caption{A correction system solver algorithm for the GLS problem}
\label{alg:gls-solver}
\begin{algorithmic}[1]
\REQUIRE  A GQR factorization of \((W,V)\) with the upper triangular matrix
\(R\in\mathbb{C}^{m\times m}\), the upper rectangular matrix \(T\in\mathbb{C}^{n\times p}\), and unitary
matrices \(Z\in\mathbb{C}^{p\times p}\), \(Q\in\mathbb{C}^{n\times n}\);
vectors \(f_{1}\in\mathbb{C}^{p}\), \(f_{2}\in\mathbb{C}^{n}\),
and \(f_{3}\in\mathbb{C}^{m}\).
\ENSURE  Vectors \(x_{1}\in\mathbb{C}^{p}\), \(x_{2}\in\mathbb{C}^{n}\),
and \(x_{3}\in\mathbb{C}^{m}\) satisfying
\[
\begin{bmatrix} I_{p} & V\herm & 0 \\ V & 0 & W \\ 0 & W\herm & 0
\end{bmatrix}\begin{bmatrix}x_{1} \\ -x_{2} \\ x_{3}\end{bmatrix}
=\begin{bmatrix}f_{1} \\ f_{2} \\ f_{3}\end{bmatrix}.
\]
\hspace{-2\algorithmicindent}
\textbf{Interface}: \((x_{1},x_{2},x_{3})=\mathtt{CorrectionSystemSolverGLS}
(R,T,Z,Q,f_{1},f_{2},f_{3})\)

\medskip
\STATE  \(u\gets Q\herm f_{2}\), i.e., \(u_{1}=(Q\herm f_{2})(1:m)\)
and \(u_{2}=(Q\herm f_{2})(m+1:n)\).
\STATE  \(w\gets Zf_{1}\), i.e., \(w_{1}=(Zf_{1})(1:p-n+m)\)
and \(w_{2}=(Zf_{1})(p-n+m+1:p)\).
\STATE  Compute \(h_{1}\) by solving the triangular
system \(R\herm h_{1}=f_{3}\).
\STATE  Compute \(g_{2}\) by solving the triangular
system \(T_{22}\herm g_{2}=u_{2}\).
\STATE  Compute \(h_{2}\) by solving the triangular
system \(T_{22}\herm h_{2}=w_{2}-g_{2}-T_{12}\herm h_{1}\).
\STATE  \(g_{1}\gets w_{1}-T_{11}\herm h_{1}\).
\STATE  \(x_{1} \gets Z\herm g\) and \(x_{2}\gets Qh\).
\STATE  Compute \(x_{3}\) by solving the triangular
system \(Rx_{3}=u_{1}-T_{11}g_{1}-T_{12}g_{2}\).
\RETURN  \((x_{1},x_{2},x_{3})\).
\end{algorithmic}
\end{algorithm}

\begin{algorithm}[!tb]
\caption{Mixed precision GLS algorithm}
\label{alg:mpgls}
\begin{algorithmic}[1]
\REQUIRE  Matrices \(W\in\mathbb{C}^{n\times m}\), \(V\in\mathbb{C}^{n\times
p}\) and vector \(d\in\mathbb{C}^{n}\) storing at precision \(\machepso\);
the maximal number of iterations \texttt{maxit} with a default value \(40\);
the tolerance parameter \texttt{tol} for the stopping criteria.
\ENSURE  Vectors \(x\in\mathbb{C}^{m}\) and \(y\in\mathbb{C}^{p}\)
stored at precision \(\machepso\) approximating
the solution of the GLS problem~\eqref{problem:gls}.

\medskip
\STATE  Compute the GQR factorization of \((W,V)\) at
precision \(\machepsf\), i.e., \(W=QR_{W}\), \(V=QTZ\).
\STATE  Compute the initial guess \(x_{0}\), \(y_{0}\) via
the GQR factorization at precision \(\machepsf\), and store
\(x\gets x_{0}\), \(y\gets y_{0}\) at precision \(\machepso\).
\STATE  Compute \(z_{0}\gets Q\herm\bar{v}\),
where \(\bar{v}=[0,v\herm]\herm\) and \(v\) is computed by solving
the triangular system \(T_{22}\herm v=(Zy)(p-n+m+1:p)\) at precision
\(\machepsf\), and store \(z\gets z_{0}\) at precision \(\machepso\).
\FOR{\(i=1,\dotsc,\mathtt{maxit}\)}
\STATE  Compute the residuals \(f_{1}\gets -y+V\herm z\),
\(f_{2}\gets d-Wx-Vy\), and \(f_{3}\gets W\herm z\) at
precision~\(\machepsr\), and store them at precision \(\machepso\).
\IF{the stopping criteria~\eqref{eq:gls-stop} are satisfied}
\RETURN  \(x\), \(y\).
\ENDIF
\STATE  Compute \(\Delta y\), \(\Delta z\), and \(\Delta x\)
at precision \(\machepss\) by function
\[
(\Delta y,\Delta z,\Delta x)=\mathtt{CorrectionSystemSolverGLS}
(R,T,Z,Q,f_{1},f_{2},f_{3})
\]
in Algorithm~\ref{alg:gls-solver}, and store them at precision \(\machepso\).
\STATE  Update and store \(x\gets x+\Delta x\), \(y\gets y+\Delta y\),
and \(z\gets z+\Delta z\) at precision \(\machepso\).
\ENDFOR
\end{algorithmic}
\end{algorithm}

Similarly to Section~\ref{subsec:mpalgo-lse}, we consider the complexity
of the algorithm by calculating the number of floating-point operations.
The cost of Algorithm~\ref{alg:mpgls} consists of
\[
t_{f}=\mathcal{O}(m^{2}n+n^{2}p)
\]
for the matrix factorization, as well as
\begin{align*}
t_{w} & =\mathcal{O}(m+n+p),\quad t_{r}=4(mn+np)+\mathcal{O}(m+n+p),\\
t_{s} & =2(pn+n^{2})+4mp+4(m^{2}+(n-m)^{2})+\mathcal{O}(m+n+p)
\end{align*}
in each refinement step, where \(t_{f}\), \(t_{w}\), \(t_{r}\), and
\(t_{s}\) represent the total number of floating-point operations at precision
\(\machepsf\), \(\machepso\), \(\machepsr\), and \(\machepss\), respectively.
If \(p\geq n\geq m\), we obtain
\[
t_{f}=\mathcal{O}(n^{2}p),\qquad t_{w}=\mathcal{O}(p),\qquad
t_{r}=\mathcal{O}(np),\qquad t_{s}=\mathcal{O}(np).
\]
The dominating term is \(t_{f}\) when the matrix dimension
\(n\) is relatively large.
\section{Accuracy analysis}
\label{sec:accuracy}
In this section, we conduct an accuracy analysis
of our classical iterative refinement algorithms.
It comprises a backward and a forward error analysis for each problem.
We state that if LSE and GLS problems are solved via the methods
in Algorithms~\ref{alg:mplse} and~\ref{alg:mpgls},
then the backward error in the correction solve is limited by
the precision in which the GRQ or GQR factorization is computed.

We first introduce some notation regarding rounding error analysis.
For simplicity, we use~\(\hat{\cdot}=\fl(\cdot)\) to represent
the computed values interchangeably without further clarification.
Let \(\gamma_{n}=n\machepso/(1-n\machepso)\) and \(\tilde{\gamma}_{n}
=cn\machepso/(1-cn\machepso)\), where \(c\) is a small constant.
Here we assume \(n\machepso\ll 1\).
The superscripts \(f\), \(s\), and \(r\) stand for the procedure computed
at precision \(\machepsf\), \(\machepss\), and \(\machepsr\), respectively.

We start by providing some assumptions and notations
used in our analysis in Section~\ref{subsec:accuracy:assump}.
Then we shall present our analysis in several parts.
In Section~\ref{subsec:accuracy:grqgqr}, we introduce
a lemma on the accuracy of the GRQ/GQR factorization.
Then in each subsection we conduct a rounding error
analysis of our mixed precision algorithm.
We summarize our results in Theorems~\ref{thm:lse} and~\ref{thm:gls}
and derive the limiting accuracy in Theorems~\ref{thm:lse-limit}
and~\ref{thm:gls-limit}, respectively.
We then present several examples to illustrate our findings.
The analysis is partly based on the work of Bj\"{o}rck in~\cite{Bjorck1967}.

\subsection{Basic assumptions}
\label{subsec:accuracy:assump}
Before conducting the analysis, we outline a few assumptions regarding matrix
and vector operations used in our classical iterative refinement algorithms.
These assumptions are aligned with those results presented in~\cite{Higham2002}.

\begin{assumption}
\label{assump:gemv}
Assume that \(A\in\mathbb{C}^{m\times n}\), \(x\in\mathbb{C}^{n}\),
and \(y\in\mathbb{C}^{m}\) so that \(y=Ax\).
For the matrix--vector multiplication, i.e., \(\mathtt{xGEMV}\) in
BLAS/LAPACK~\cite{lapack1999}, the computed result \(\hat{y}=\fl(Ax)\) satisfies
\[
\hat{y}=(A+\Delta A)x,\qquad\abs{\Delta A}\leq\epsgemv(n)\abs{A},
\]
where \(\epsgemv(n)\) denotes a function of \(\machepso\) and \(n\).
\end{assumption}

\begin{assumption}
\label{assump:qr}
Assume that \(A\in\mathbb{C}^{m\times n}\).
For the QR factorization \(A=QR\), i.e., \(\mathtt{xGEQRF}\)
in LAPACK~\cite{lapack1999}, with unitary \(Q\in\mathbb{C}^{m\times m}\)
and upper rectangular \(R\in\mathbb{C}^{m\times n}\), there exists
a unitary matrix \(Q_{0}\in\mathbb{C}^{m\times m}\) such that
\[
A+\Delta A=Q_{0}\hat{R},\qquad\norm{\Delta a_{j}}\leq\epsqr(m,n)\norm{a_{j}},
\]
where \(\hat{R}\) is the computed triangular factor, \(A=[a_{1},\dotsc,a_{n}]\),
and \(\epsqr(m,n)\) denotes a function of \(\machepso\), \(m\) and \(n\).
\end{assumption}

\begin{assumption}
\label{assump:trsv}
Assume that \(T\in\mathbb{C}^{n\times n}\) is a triangular matrix,
and \(b\in\mathbb{C}^{n}\).
For the triangular system solver \(Tx=b\), i.e., \(\mathtt{xTRSV}\) in
BLAS/LAPACK~\cite{lapack1999}, the computed result \(\hat{x}\) satisfies
\[
(T+\Delta T)\hat{x}=b,\qquad\abs{\Delta T}\leq\epstrsv(n)\abs{T},
\]
where \(\epstrsv(n)\) denotes a function of \(\machepso\) and \(n\).
\end{assumption}

\begin{assumption}
\label{assump:unmqr}
Assume that \(Q\in\mathbb{C}^{m\times m}\) is a unitary matrix,
and \(y\), \(z\in\mathbb{C}^{m}\) so that \(y=Qz\).
For applying \(Q\) to \(z\), i.e.,
\(\mathtt{xORMQR}\slash{}\mathtt{xUNMQR}\) in LAPACK~\cite{lapack1999},
there exist a backward error \(\Delta z\) and a forward error
\(\Delta y\) such that the computed result \(\hat{y}=\fl(Qz)\) satisfies
\begin{align*}
\hat{y}=Q(z+\Delta z),\qquad & \norm{\Delta z}\leq\epsunmqr(m)\norm{z},\\
\hat{y}=Qz+\Delta y,\qquad & \norm{\Delta y}\leq\epsunmqr(m)\norm{\hat{y}},
\end{align*}
where \(\epsunmqr(m)\) denotes a function of \(\machepso\) and \(m\).
\end{assumption}

In Table~\ref{tab:notation}, we present several established
results related to these \(\epsilon_{*}(\cdot)\)'s.
Nonetheless, stricter bounds are available, prompting the use
of the more general notation \(\epsilon_{*}(\cdot)\).

\begin{table}[!tb]
\centering
\caption{Notations of rounding error analysis.
Here \(A\in\mathbb{C}^{m\times n}\); \(Q\), \(Q_{0}\in\mathbb{C}^{m\times m}\)
are unitary matrices; \(T\in\mathbb{C}^{n\times n}\) is a triangular
matrix; \(R\in\mathbb{C}^{m\times n}\) is an upper rectangular matrix;
\(b\), \(x\in\mathbb{C}^{n}\), and \(y\), \(z\in\mathbb{C}^{m}\) are vectors.}
\label{tab:notation}
\begin{spacing}{1.3}
\medskip
\begin{tabular}{ccccc}
\hline Notation & LAPACK Subroutine & Problem
& Upper bound on \(\epsilon\) & Sources \\ \hline
\(\epsgemv(n)\)  & \texttt{xGEMV} & \(y=Ax\)
& \(\gamma_{n}\) & \cite[Section 3.5]{Higham2002}\\
\(\epsqr(m,n)\)  & \texttt{xGEQRF} & \(A=QR\)
& \(\tilde{\gamma}_{mn}\)    & \cite[Section 19.3]{Higham2002}\\
\(\epstrsv(n)\)  & \texttt{xTRSV} & \(Tx=b\)
& \(\gamma_{n}\) & \cite[Section 8.1]{Higham2002}\\
\(\epsunmqr(m)\) & \texttt{xUNMQR} & \(y=Qz\)
& \(\tilde{\gamma}_{m^{2}}\) & \cite[Section 19.3]{Higham2002}\\ \hline
\end{tabular}
\end{spacing}
\end{table}

Employing pre-defined notations from Table~\ref{tab:notation}
instead of existing error bounds like \(\gamma_{n}\)
is beneficial to the discussion of rounding errors.
On one hand, it helps us to trace back and locate
where the dominant rounding error terms come from and
provides better insight for mixed precision algorithms.
On the other hand, bounds like \(\gamma_{n}\) are often an overestimate.
It may fail to provide useful information if the matrix
dimension \(n\) is large and the working precision is too low.
Therefore it is recommended to use such pre-defined
notations in future rounding error analysis.

\subsection{Analysis for the GRQ\slash{}GQR factorization}
\label{subsec:accuracy:grqgqr}
A detailed review on the GRQ factorization can be found in~\cite{ABD1992}.
In general, it comprises two matrix factorizations: an RQ
factorization \(B=[0,R]Q\) and then a QR factorization \(AQ\herm=ZT\).
The LAPACK subroutine \texttt{xGGRQF} is also implemented in this manner.

\begin{lemma}
\label{lemma:grq}
Let \(A\in\mathbb{C}^{m\times n}\) and \(B\in\mathbb{C}^{p\times n}\).
Suppose \((\hat{R},\hat{T})\) are the computed triangular matrices from
the GRQ factorization of \((B,A)\) via Householder transformations.
Then there exist matrices \(\Delta E_{1}\), \(\Delta E_{2}\),
and unitary matrices \(Q\), \(Z\) such that
\[
B+\Delta E_{1}=[0,\hat{R}]Q,\qquad A+\Delta E_{2}=Z\hat{T}Q,
\]
where \(\normfro{\Delta E_{1}}\leq\epsqr(n,p)\normfro{B}\) and
\[
\normfro{\Delta E_{2}}\leq\big(\epsunmqr(n)+\epsqr(m,n)
+\epsunmqr(n)\epsqr(m,n)\big)\normfro{A}.
\]
\end{lemma}

\begin{proof}
The GRQ process first conducts an RQ factorization \(B=[0,R]Q\).
According to Assumption~\ref{assump:qr}, there exists a triangular matrix
\(\hat{R}\) such that the backward error \(\Delta E_{1}=[0,\hat{R}]Q-B\)
satisfies \(\normfro{\Delta E_{1}}\leq\epsqr(n,p)\normfro{B}\).
Since \(Q\) is unitary, we have
\begin{equation}
\label{eq:lem-GRQ-AQherm}
\bignormfro{\fl(AQ\herm)-AQ\herm}\leq\epsunmqr(n)\normfro{A}.
\end{equation}
Thus
\[
\bignormfro{\fl(AQ\herm)}\leq\normfro{AQ\herm}
+\bignormfro{\fl(AQ\herm)-AQ\herm}\leq(1+\epsunmqr(n))\normfro{A}.
\]
For the QR factorization of \(AQ\herm=ZT\), an error bound is presented as
\begin{equation}
\label{eq:lem-qr}
\bignormfro{Z\hat{T}-\fl(AQ\herm)}\leq\epsqr(m,n)
\normfro{\fl(AQ\herm)}\leq\epsqr(m,n)(1+\epsunmqr(n))\normfro{A}.
\end{equation}
Combining~\eqref{eq:lem-GRQ-AQherm} and~\eqref{eq:lem-qr}, we derive
that \(\Delta E_{2}=Z\hat{T}Q-A=(Z\hat{T}-AQ\herm)Q\) satisfies
\[
\normfro{\Delta E_{2}}=\normfro{Z\hat{T}-AQ\herm}\leq\big(\epsunmqr(n)
+\epsqr(m,n)+\epsunmqr(n)\epsqr(m,n)\big)\normfro{A}.\qedhere
\]
\end{proof}

Likewise, the rounding error analysis of the GQR factorization
is stated as the following Lemma~\ref{lemma:gqr}.
The proof is essentially the same and thus omitted.

\begin{lemma}
\label{lemma:gqr}
Let \(W\in\mathbb{C}^{n\times m}\) and \(V\in\mathbb{C}^{n\times p}\).
Suppose \((\hat{R},\hat{T})\) are the computed triangular matrices from
the GQR factorization of \((W,V)\) via Householder transformations.
Then there exist matrices \(\Delta F_{1}\), \(\Delta F_{2}\),
and unitary matrices \(Q\), \(Z\) such that
\[
W+\Delta F_{1}=Q\begin{bmatrix}\hat{R} \\ 0\end{bmatrix},\qquad V+\Delta F_{2}=Q\hat{T}Z,
\]
where \(\normfro{\Delta F_{1}}\leq\epsqr(n,m)\normfro{W}\) and
\[
\normfro{\Delta F_{2}}\leq\bigl(\epsunmqr(n)
+\epsqr(p,n)+\epsunmqr(n)\epsqr(p,n)\bigr)\normfro{V}.
\]
\end{lemma}

\begin{subsection}{Analysis for the LSE problem}
\label{subsec:accuracy-lse}
We now state the main theorem of our analysis.

\begin{theorem}
\label{thm:lse}
Let \(A\in\mathbb{C}^{m\times n}\) and \(B\in\mathbb{C}^{p\times n}\).
Suppose the correction system~\eqref{eq:correction-system-lse}
is solved via Algorithm~\ref{alg:lse-solver}.
Then there exist \(\Delta A_{1}\), \(\Delta A_{2}\), \(\Delta B_{1}\), \(\Delta
B_{2}\), and \(\Delta t_{1}\), \(\Delta t_{2}\), \(\Delta t_{3}\) such that
\[
\begin{bmatrix}I_{m} & 0 & A+\Delta A_{1} \\ 0 & 0 & B+\Delta B_{1}\\
(A+\Delta A_{2})\herm & (B+\Delta B_{2})\herm & 0\end{bmatrix}\begin{bmatrix}
\fl(\Delta r) \\ -\fl(\Delta v) \\ \fl(\Delta x)\end{bmatrix}=\begin{bmatrix}
f_{1}+\Delta t_{1} \\ f_{2}+\Delta t_{2} \\ f_{3}+\Delta t_{3}\end{bmatrix},
\]
and the following bounds are satisfied:
\begin{align*}
\normfro{\Delta A_{i}}\leq{} &
(\eta_{0}+\eta_{i}+\eta_{0}\eta_{i})\normfro{A},\quad i=1,2,\\
\normfro{\Delta B_{i}}\leq{} & \bigl(\epsqr^{f}(n,p)+5\epstrsv^{s}(p)
+5\epsqr^{f}(n,p)\epstrsv^{s}(p)\bigr)\normfro{B},\quad i=1,2,\\
\norm{\Delta t_{1}}\leq{} & \epsunmqr^{s}(m)\norm{f_{1}}
+3\epsunmqr^{s}(m)\norm{\fl(\Delta r)}+(1+\eta_{0})(1+\eta_{1})\,
\epsunmqr^{s}(n)\normfro{A}\norm{\fl(\Delta x)},\\
\norm{\Delta t_{2}}\leq{} & (1+\epsqr^{f}(n,p))(1+\epstrsv^{s}(p))\,
\epsunmqr^{s}(n)\normfro{B}\norm{\fl(\Delta x)},\\
\norm{\Delta t_{3}}\leq{} & \epsunmqr^{s}(n)\norm{f_{3}}
+3(1+\eta_{0})(1+\eta_{2})\,\epsunmqr^{s}(m)\normfro{A}\norm{\fl(\Delta r)},
\end{align*}
where
\begin{align*}
\eta_{0} & =\epsunmqr^{f}(n)+\epsqr^{f}(m,n)+\epsunmqr^{f}(n)\epsqr^{f}(m,n),\\
\eta_{1} & =\max\bigl\{\epsgemv^{s}(p),4\epstrsv^{s}(n-p)\bigr\},\\
\eta_{2} & =\max\bigl\{\epstrsv^{s}(n-p),
\epsgemv^{s}(n-p),\epsgemv^{s}(m-n+p)\bigr\}
\end{align*}
are small constants of \(m\), \(n\), \(p\), and the unit roundoffs.
\end{theorem}

\begin{proof}
The proof of this theorem is quite long and tedious.
We present it in Appendix~\ref{sec:appendix-proof}.
\end{proof}

\begin{remark}
If we use the estimates in Table~\ref{tab:notation},
the upper bounds of errors are simplified as
\begin{align*}
& \normfro{\Delta A_{i}}=\mathcal{O}(\machepsf)\,\normfro{A},\qquad
\normfro{\Delta B_{i}}=\mathcal{O}(\machepsf)\,\normfro{B},\quad i=1,2,\\
& \norm{\Delta t_{1}}=\mathcal{O}(\machepss)\,\bigl(\norm{f_{1}}+
\norm{\fl(\Delta r)}+\normfro{A}\norm{\fl(\Delta x)}\bigr),\\
& \norm{\Delta t_{2}}=\mathcal{O}(\machepss)\,\normfro{B}\norm{\fl(\Delta x)},\\
& \norm{\Delta t_{3}}=\mathcal{O}(\machepss)\,\bigl(\norm{f_{3}}
+\normfro{A}\norm{\fl(\Delta r)}\bigr).
\end{align*}
\end{remark}

Theorem~\ref{thm:lse} essentially says that the backward error in the correction
solve is limited by the precision in which the GRQ factorization is computed.
As in Section~\ref{sec:algorithm}, for the LSE problem, suppose
\begin{equation}
\label{eq:lse-limit-notation}
\tilde{F}=\begin{bmatrix}\alpha I_{m} & 0 & A \\ 0 & 0 & \beta B \\
A\herm & \beta B\herm & 0\end{bmatrix},\quad
\tilde{u}=\begin{bmatrix}r \\ -v \\ x\end{bmatrix},\quad
\Delta\tilde{u}=\begin{bmatrix}\alpha^{-1/2}\Delta r \\
-\alpha^{-1/2}\beta^{-1}\Delta v \\ \alpha^{1/2}\Delta x\end{bmatrix},\quad
\tilde{s}=\begin{bmatrix}\alpha^{1/2}f_{1} \\
\alpha^{1/2}\beta f_{2} \\ \alpha^{-1/2}f_{3}\end{bmatrix},
\end{equation}
where \(\alpha\) and \(\beta\) are positive scaling parameters.
Then system~\eqref{eq:correction-system-lse} is
equivalent to \(\tilde{F}\Delta\tilde{u}=\tilde{s}\),
and the refinement stage of Algorithm~\ref{alg:mplse} can be regarded
as a solver of the linear system \(\tilde{F}\Delta\tilde{u}=\tilde{s}\).
We have the following result.

\begin{theorem}
\label{thm:lse-limit}
For matrices \(A\) and \(B\) satisfying \(\machepsf\kappa_{\infty}
(\tilde{F})\lesssim 1\) and any scaling parameter \(\alpha\),~\(\beta>0\),
the limiting relative forward error is derived as
\begin{equation}
\label{eq:lse-limit}
\frac{\norminf{\tilde{u}-\fl(\tilde{u})}}{\norminf{\tilde{u}}}
\leq4(m+n+p+1)\cond(\tilde{F},\tilde{u})\machepsr+\machepso,
\end{equation}
where \(\cond(\tilde{F},\tilde{u})=\bignorminf{\abs{\tilde{F}^{-1}}
\cdot\abs{\tilde{F}}\cdot\abs{\tilde{u}}}/\norminf{\tilde{u}}\).
\end{theorem}

\begin{proof}
The proof is similar to~\cite[Section 2]{CHP2020} and thus we omit it.
\end{proof}

\begin{remark}
\label{remark:lse}
A key benefit of iterative refinement is that the accuracy of
Algorithm~\ref{alg:mplse} is \emph{independent} of scaling.
In other words, users are not required to select appropriate
\(\alpha\) and \(\beta\) so as to balance the condition number
of \(\tilde{F}\) from an algorithmic perspective.
At the same time, the result holds for the optimal choice of
\(\alpha\) and \(\beta\) that minimizes the condition number
of \(\tilde{F}\) from a theoretical standpoint.
\end{remark}
\end{subsection}

\begin{subsection}{Analysis for the GLS problem}
\label{subsec:accuracy-gls}
The accuracy analysis for the GLS problem is similar to
the arguments in Section~\ref{subsec:accuracy-lse}.
The proofs of the theorems are omitted.

\begin{theorem}
\label{thm:gls}
Let \(W\in\mathbb{C}^{n\times m}\) and \(V\in\mathbb{C}^{n\times p}\).
Suppose the correction system~\eqref{eq:correction-system-gls}
is solved via Algorithm~\ref{alg:gls-solver}.
Then there exist \(\Delta W_{1}\), \(\Delta W_{2}\), \(\Delta V_{1}\), \(\Delta
V_{2}\), and \(\Delta t_{1}\), \(\Delta t_{2}\), \(\Delta t_{3}\) such that
\[
\begin{bmatrix}I_{p} & (V+\Delta V_{1})\herm & 0 \\ V+\Delta V_{2} & 0
& W+\Delta W_{2} \\ 0 & (W+\Delta W_{1})\herm & 0\end{bmatrix}\begin{bmatrix}
\fl(\Delta y) \\ -\fl(\Delta z) \\ \fl(\Delta x)\end{bmatrix}=\begin{bmatrix}
f_{1}+\Delta t_{1} \\ f_{2}+\Delta t_{2} \\ f_{3}+\Delta t_{3}\end{bmatrix},
\]
and the following bounds are satisfied:
\begin{align*}
\normfro{\Delta V_{i}}\leq{} &
(\eta_{0}+\eta_{i}+\eta_{0}\eta_{i})\normfro{V},\quad i=1,2,\\
\normfro{\Delta W_{i}}\leq{} & \bigl(\epsqr^{f}(n,m)+5\epstrsv^{s}(m)
+5\epsqr^{f}(n,m)\epstrsv^{s}(m)\bigr)\normfro{W},\quad i=1,2,\\
\norm{\Delta t_{1}}\leq{} & \epsunmqr^{s}(p)\norm{f_{1}}+3\epsunmqr^{r}(p)
\norm{\fl(\Delta y)}+(1+\eta_{0})(1+\eta_{1})\,
\epsunmqr^{r}(n)\normfro{V}\norm{\fl(\Delta z)},\\
\norm{\Delta t_{2}}\leq{} & \epsunmqr^{s}(n)\norm{f_{2}}+3(1+\eta_{0})
(1+\eta_{2})\,\epsunmqr^{r}(p)\normfro{V}\norm{\fl(\Delta y)},\\
\norm{\Delta t_{3}}\leq{} & (1+\epsqr^{f}(n,m))(1+\epstrsv^{s}(m))\,
\epsunmqr^{r}(n)\normfro{W}\norm{\fl(\Delta z)},
\end{align*}
where
\begin{align*}
\eta_{0} & =\epsunmqr^{f}(n)+\epsqr^{f}(p,n)+\epsunmqr^{f}(n)\epsqr^{f}(p,n),\\
\eta_{1} & =\max\big\{\epsgemv^{s}(m),4\epstrsv^{s}(n-m)\big\},\\
\eta_{2} & =\max\big\{\epstrsv^{s}(n-m),
\epsgemv^{s}(n-m),\epsgemv^{s}(p-n+m)\big\}
\end{align*}
are small constants of \(m\), \(n\), \(p\), and the unit roundoffs.
\end{theorem}

\begin{remark}
If we use the estimates in Table~\ref{tab:notation},
the upper bounds of errors are simplified as
\begin{align*}
& \normfro{\Delta V_{i}}=\mathcal{O}(\machepsf)\,\normfro{V},\qquad
\normfro{\Delta W_{i}}=\mathcal{O}(\machepsf)\,\normfro{W},\quad i=1,2,\\
& \norm{\Delta t_{1}}=\mathcal{O}(\machepss)\,\bigl(\norm{f_{1}}
+\norm{\fl(\Delta y)}+\normfro{V}\norm{\fl(\Delta z)}\bigr),\\
& \norm{\Delta t_{2}}=\mathcal{O}(\machepss)\,\bigl(\norm{f_{2}}
+\normfro{V}\norm{\fl(\Delta y)}\bigr),\\
& \norm{\Delta t_{3}}=\mathcal{O}(\machepss)\,\normfro{W}\norm{\fl(\Delta z)}.
\end{align*}
\end{remark}

As in Section~\ref{sec:algorithm}, suppose
\begin{equation}
\label{eq:gls-limit-notation}
\tilde{F}=\begin{bmatrix}\alpha I_{p} & V\herm & 0
\\ V & 0 & \beta W \\ 0 & \beta W\herm & 0\end{bmatrix},\quad
\tilde{u}=\begin{bmatrix}y \\ -z \\ x\end{bmatrix},\quad
\Delta\tilde{u}=\begin{bmatrix}\alpha^{-1/2}\Delta y \\
-\alpha^{1/2}\Delta z \\ \alpha^{-1/2}\beta^{-1}\Delta x\end{bmatrix},\quad
\tilde{s}=\begin{bmatrix}\alpha^{1/2}f_{1} \\
\alpha^{-1/2}f_{2} \\ \alpha^{1/2}\beta f_{3}\end{bmatrix},
\end{equation}
where \(\alpha\) and \(\beta\) are positive scaling parameters.
Then \(\tilde{F}\Delta\tilde{u}=\tilde{s}\).
Similarly to Theorem~\ref{thm:lse-limit}, the limiting
accuracy of the GLS problem is presented as follows.

\begin{theorem}
\label{thm:gls-limit}
For matrices \(W\) and \(V\) that satisfy \(\machepsf\kappa_{\infty}
(\tilde{F})\lesssim 1\), the limiting relative forward error is derived as
\begin{equation}
\label{eq:gls-limit}
\frac{\norminf{\tilde{u}-\fl(\tilde{u})}}{\norminf{\tilde{u}}}
\leq 4(m+n+p+1)\cond(\tilde{F},\tilde{u})\machepsr+\machepso.
\end{equation}
\end{theorem}

Theorems~\ref{thm:lse-limit} and~\ref{thm:gls-limit} imply that
\(\kappa_{\infty}(\tilde{F})\) must be bounded by \(\mathcal{O}
(\machepsf^{-1})\) in order for the algorithm to converge,
and precision \(\machepso\) together with the product
\(\cond(\tilde{F},\tilde{u})\machepsr\) determine the attainable accuracy.
A specific example of that can be found in Table~\ref{tab:qr} as
in~\cite{CH2018} by employing IEEE half/single/double precisions.

\begin{table}[!tb]
\begin{spacing}{1.3}
\caption{Comparison of results for classical iterative refinement.}
\label{tab:qr}
\medskip\centering
\begin{tabular}{ccccc}\hline
\(\machepsf\) & \(\machepso\) & \(\machepsr\)
& Upper bound on \(\kappa_{\infty}(\tilde{F})\) & Forward error\\ \hline
single & single & double & \(\approx10^{8}\) & single\\
single & double & double & \(\approx10^{8}\)
& \(\cond(\tilde{F},\tilde{u})\cdot10^{-16}\) \\
single & double &  quad  & \(\approx10^{8}\) & double \\
 half  & single & double & \(\approx10^{4}\) & single\\
 half  & double & double & \(\approx10^{4}\)
& \(\cond(\tilde{F},\tilde{u})\cdot10^{-16}\) \\
 half  & double &  quad  & \(\approx10^{4}\) & double \\ \hline
\end{tabular}
\end{spacing}
\end{table}

Another conclusion of Theorems~\ref{thm:lse-limit} and~\ref{thm:gls-limit}
is that, if we use a higher (than working) precision to compute the residual
in the iterative refinement, i.e., \(\machepsr\ll\machepso\),
then \(\cond(\tilde{F},\tilde{u})\machepsr\) may no longer
be the leading error term in~\eqref{eq:lse-limit} and~\eqref{eq:gls-limit}.
The forward error is thus reduced from \(\cond(\tilde{F},\tilde{u})
\cdot10^{-16} \) to double precision, as illustrated
by certain scenarios in Table~\ref{tab:qr}.
In other words, the mixed precision algorithm
improves the accuracy of the solution.
\end{subsection}
\section{GMRES-based iterative refinement}
\label{sec:gmres}
As shown in Section~\ref{sec:accuracy}, the mixed precision classical
iterative refinement algorithm has the obvious drawback that it only
works for the case \(\machepsf\kappa_{\infty}(\tilde{F})\lesssim 1\).
To overcome this problem, we use preconditioned GMRES
algorithms to solve the saddle-point system.
Since the GRQ or GQR factorization is already obtained, we may
construct the preconditioners using these computed factors.
For a brief introduction to the GMRES algorithm, see~\cite{SS1986}.

Let \(\tilde{F}\) be the augmented matrix with
scaling parameters \(\alpha\) and \(\beta\), i.e.,
\[
\tilde{F}=\begin{bmatrix}\alpha I_{m} & 0 & A \\ 0 & 0 & \beta B \\
A\herm & \beta B\herm & 0\end{bmatrix}\quad\text{(LSE)}\qquad\text{or}\qquad
\tilde{F}=\begin{bmatrix}\alpha I_{p} & V\herm & 0 \\ V & 0
& \beta W \\ 0 & \beta W\herm & 0\end{bmatrix}\quad\text{(GLS)}.
\]
Here the parameters \(\alpha\) and \(\beta\) are designed to reduce
the condition number of \(\tilde{F}\) and to balance the norms of matrices
\(A\) and \(B\) for the LSE problem (\(V\) and \(W\) for the GLS problem).
Determining the optimal values of \(\alpha\) and \(\beta\)
for the LSE and GLS problems is theoretically challenging,
as the spectrum of \(\tilde{F}\) lacks clarity.
For LSE problems, during numerical experiments, we set
\(\beta\approx\norm{A}/\norm{B}\) to achieve norm balance,
and our choice of \(\alpha\) is derived from the particular scenario \(B=0\),
where the LSE problem simplifies to the standard least squares problem.
Under the assumption \(B=0\), the optimal choice of \(\alpha\) is
\(\alpha=2^{-1/2}\sigma_{\min}(A)\)~\cite{Bjorck1967,CD2025}.
Since computing \(\sigma_{\min}(A)\) is costly,
we instead use \(\alpha=\norm{r_{0}}\),
a choice inspired by~\cite{ADR1989} in our numerical experiments.
Similarly, we set \(\alpha=\norm{y_{0}}\) and
\(\beta\approx\norm{V}/\norm{W}\) for GLS problems.

These two augmented systems with \(\tilde{F}\) can be
regarded as symmetric saddle-point systems.
In~\cite{BGL2005,Rozloznik2018}, both left and block diagonal
preconditioners for saddle-point systems are studied.
Here we propose an improved version of preconditioners
for addressing LSE and GLS problems.

\subsection{Left preconditioners}
\label{subsec:gmres-left}
A natural choice of the left preconditioner for the LSE problem,
inspired by directly calculating the inverse of \(\tilde{F}\), is given as
\begin{equation}
\label{eq:lse-leftpre}
M=\begin{bmatrix}\alpha^{-1}I_{m} & -\alpha^{-1}\beta^{-1}AB\pinv & 0 \\
-\alpha^{-1}\beta^{-1}(AB\pinv)\herm & \alpha^{-1}\beta^{-2}
(B\pinv)\herm A\herm AB\pinv & \beta^{-1}(B\pinv)\herm \\
0 & \beta^{-1}B\pinv & 0\end{bmatrix}.
\end{equation}

Similarly, a left preconditioner for the GLS problem is presented as
\begin{equation}
\label{eq:gls-leftpre}
M=\begin{bmatrix}\alpha^{-1}I_{p} & 0 & -\alpha^{-1}\beta^{-1}V\herm
(W\pinv)\herm \\ 0 & 0 & \beta^{-1}(W\pinv)\herm \\
-\alpha^{-1}\beta^{-1}W\pinv V & \beta^{-1}W\pinv
& \alpha^{-1}\beta^{-2}W\pinv VV\herm(W\pinv)\herm\end{bmatrix}.
\end{equation}
For both LSE and GLS problems,
applying the preconditioner \(M\) to a vector is conducted using
the procedure in Section~\ref{sec:algorithm} in practice.

The convergence of the GMRES algorithm with these left
preconditioners follows from the arguments in~\cite{CHP2020}.
However, as shown in Section~\ref{sec:experiments}, these specific left
preconditioners do not work well in the mixed precision GMRES-based iterative
refinement algorithm and sometimes take an awfully long time to converge
as each iteration is numerically very expensive.
This motivates us to develop other preconditioners
in order to speed up the process.

\subsection{Block-diagonal split preconditioner for the LSE problem}
\label{subsec:gmres-lse}
If \(m\geq n\), \(\hat{T}\) can be partitioned into two parts
\(\hat{T}\herm=[\hat{T}_{1}\herm,0]\) with \(\hat{T}_{1}
\in\mathbb{C}^{n\times n}\) being upper triangular.
Then a block-diagonal split preconditioner for the LSE problem is presented as
\begin{equation}
\label{eq:lse-bdpre-mlarge}
\begin{split}
M_{l} & =\diag\big\{\alpha^{-1/2}I_{m},\alpha^{-1/2}\beta^{-1}
S\hat{R}^{-1},\alpha^{1/2}\hat{T}_{1}\iherm Q\big\},\\
M_{r} & =\diag\big\{\alpha^{-1/2}I_{m},\alpha^{-1/2}\beta^{-1}
\hat{R}\iherm S\herm,\alpha^{1/2}Q\herm\hat{T}_{1}^{-1}\big\},
\end{split}
\end{equation}
where \(S=\hat{T}_{1}(n-p+1:n,n-p+1:n)\) is a \(p\)-by-\(p\) matrix.

We aim to derive an upper bound for the condition
number \(\kappa_{\infty}(M_{l}\tilde{F}M_{r})\).
To better reflect the actual case we consider \(\tilde{F}\)
with rounding error during the process, i.e.,
\(B+\Delta E_{1}=[0,\hat{R}]Q\), and \(A+\Delta E_{2}=Z\hat{T}Q\).
In exact arithmetic one can simply take \(\Delta E_{i}=0\) in the analysis.
Let
\[
X=\begin{bmatrix}I_{m} & 0 & Z_{1} \\ 0 & 0 & \bar{I}_{p} \\
Z_{1}\herm & \bar{I}_{p}\herm & 0\end{bmatrix},\qquad\Delta E=M_{l}
\begin{bmatrix}0 & 0 &\Delta E_{2} \\ 0 & 0 & \beta\Delta E_{1} \\
\Delta E_{2}\herm & \beta\Delta E_{1}\herm & 0\end{bmatrix}M_{r},
\]
where \(\bar{I}_{p}=[0,I_{p}]\in\mathbb{C}^{p\times n}\),
\(Z=[Z_{1},Z_{2}]\) with \(Z_{1}\in\mathbb{C}^{m\times n}\).
One can verify that
\[
M_{l}\tilde{F}M_{r}=X-\Delta E.
\]
By Lemma~\ref{lemma:grq}, \(\norm{\Delta E_{1}}=\mathcal{O}(\machepsf)\norm{B}\)
and \(\norm{\Delta E_{2}}=\mathcal{O}(\machepsf)\norm{A}\).
Note that from~\eqref{eq:lse-bdpre-mlarge}, we have
\[
\norm{\Delta E}=\normbig{\begin{bmatrix}0 & 0 &
\Delta E_{2}Q\herm\hat{T}_{1}^{-1} \\
0 & 0 & S\hat{R}^{-1}\Delta E_{1}Q\herm\hat{T}_{1}^{-1} \\
\hat{T}_{1}\iherm Q\Delta E_{2}\herm & \hat{T}_{1}\iherm
Q\Delta E_{1}\herm\hat{R}\iherm S\herm & 0\end{bmatrix}}
\leq\mathcal{O}(\machepsf)\kappa_{2}(A)\kappa_{2}(B).
\]

The next step is then to derive the eigenvalues of \(X\).
A similar analysis can be found in~\cite[Section 10]{BGL2005}.
By examining the specific structure of the preconditioned matrix,
we verify that \(X\) satisfies
\[
(X^{3}-X^{2}-2X+I)(X^{2}-X-I)(X-I)=0.
\]
Therefore the spectrum of \(X\) falls into the set of
\(\Lambda\subset\{1,(1\pm\sqrt{5})/2,\lambda_{1},\lambda_{2},\lambda_{3}\}\),
where \(\lambda_{1}\approx-1.2470\), \(\lambda_{2}\approx0.4450\),
and \(\lambda_{3}\approx1.8019\) are the roots of
\(\lambda^{3}-\lambda^{2}-2\lambda+1=0\).
Since \(X\) is Hermitian, the singular values of \(X\)
equal to the absolute values of the eigenvalues.
Then it can be shown that
\[
\sigma_{\max}(X)=\lambda_{3},\qquad\sigma_{\min}(X)=\lambda_{2}.
\]
A detailed proof of this property is presented in Lemma~\ref{lemma:appendix}.
By~\cite[Corollary 2.4.4]{GV2013}, as long as \(2\norm{\Delta E}<\lambda_{2}\),
\begin{align*}
\norm{M_{l}\tilde{F}M_{r}} & =\sigma_{\max}(X-\Delta E)\leq\lambda_{3}
+\norm{\Delta E},\\ \norm{(M_{l}\tilde{F}M_{r})^{-1}} &
=\frac{1}{\sigma_{\min}(X-\Delta E)}\leq\frac{1}{\lambda_{2}-\norm{\Delta E}}.
\end{align*}
Let \(C=1+2\cdot\lambda_{3}/\lambda_{2}\approx9.0984\). Therefore
\[
\kappa_{2}(M_{l}\tilde{F}M_{r})=\norm{M_{l}\tilde{F}M_{r}}\cdot
\norm{(M_{l}\tilde{F}M_{r})^{-1}}\leq\frac{\lambda_{3}+\norm{\Delta E}}
{\lambda_{2}-\norm{\Delta E}}\leq1+2\cdot\frac{\lambda_{3}}{\lambda_{2}}=C.
\]
By the norm inequalities~\cite[Table 6.2]{Higham2002},
\[
\kappa_{\infty}(M_{l}\tilde{F}M_{r})\leq(m+n+p)
\kappa_{2}(M_{l}\tilde{F}M_{r})\leq C\cdot(m+n+p).
\]

If \(n>m\), we partition \(\hat{T}\) by \(\hat{T}=[\hat{T}_{1},\hat{T}_{2}]\)
with \(\hat{T}_{1}\in\mathbb{C}^{m\times m}\) being upper triangular
and \(\hat{T}_{2}\in\mathbb{C}^{m\times(n-m)}\).
Then a block-diagonal split preconditioner is given as
\begin{equation}
\label{eq:lse-bdpre-nlarge}
\begin{split}
M_{l} & =\diag\big\{\alpha^{-1/2}I_{m},\alpha^{-1/2}\beta^{-1}
Y\hat{R}^{-1},\alpha^{1/2}U\iherm Q\big\},\\
M_{r} & =\diag\big\{\alpha^{-1/2}I_{m},\alpha^{-1/2}\beta^{-1}
\hat{R}\iherm Y\herm,\alpha^{1/2}Q\herm U^{-1}\big\},
\end{split}
\end{equation}
where \(U\in\mathbb{C}^{n\times n}\) and
\(Y\in\mathbb{C}^{p\times p}\) are defined as
\[
U=\begin{bmatrix}\hat{T}_{1} & \hat{T}_{2} \\ 0 & I_{n-m}\end{bmatrix},
\qquad Y=\begin{bmatrix}\hat{T}_{1}(n-p+1:m,n-p+1:m)
& \hat{T}_{2}(n-p+1:m,:) \\ 0 & I_{n-m}\end{bmatrix}.
\]
One can verify that
\[
M_{l}\tilde{F}M_{r}=\begin{bmatrix}I_{m} & 0 & \bar{Z} \\ 0 & 0 & \bar{I}_{p} \\
\bar{Z}\herm & \bar{I}_{p}\herm & 0\end{bmatrix}-\Delta E,
\]
where \(\bar{Z}=[Z,0]\in\mathbb{C}^{m\times n}\).
Similarly, the condition number of the preconditioned matrix satisfies
\[
\kappa_{\infty}(M_{l}\tilde{F}M_{r})\leq(m+n+p)\cdot\frac{\lambda_{3}
+\norm{\Delta E}}{\lambda_{2}-\norm{\Delta E}}\leq C\cdot(m+n+p).
\]

Although it does not necessarily ensure a fast convergence rate theoretically,
the condition number of \(M_{l}\tilde{F}M_{r}\) has strong implications on the
magnitude of the backward and forward errors for the GMRES-based algorithms.
For references, a comprehensive bound for the backward and forward errors
in the FGMRES algorithm is provided in~\cite[Section 2]{CD2024}.
In practice, the block-diagonal split two-sided preconditioner,
due to its smaller condition number and reduced communication cost,
shows superior performance in numerical experiments, and therefore
emerges as a popular choice for GMRES-based iterative refinement.

\subsection{Block-diagonal split preconditioner for the GLS problem}
\label{subsec:gmres-gls}
If \(n\leq p\), \(\hat{T}\) can be partitioned into
two parts \(\hat{T}=[0,\hat{T}_{2}]\) with \(\hat{T}_{2}
\in\mathbb{C}^{n\times n}\) being upper triangular.
Then a two-sided preconditioner for the GLS problem is given as
\begin{equation}
\label{eq:gls-bdpre-plarge}
\begin{split}
M_{l} & =\diag\big\{\alpha^{-1/2}I_{p},\alpha^{1/2}\hat{T}_{2}^{-1}Q\herm,
\alpha^{-1/2}\beta^{-1}S\herm\hat{R}\iherm\big\},\\
M_{r} & =\diag\big\{\alpha^{-1/2}I_{p},\alpha^{1/2}Q\hat{T}_{2}\iherm,
\alpha^{-1/2}\beta^{-1}\hat{R}^{-1}S\big\},
\end{split}
\end{equation}
where \(S=\hat{T}_{2}(1:m,1:m)\) is an \(m\)-by-\(m\) matrix.
One can verify that
\begin{equation}
\label{eq:gls-case}
M_{l}\tilde{F}M_{r}\approx\begin{bmatrix}I_{p} & Z_{2}\herm & 0 \\
Z_{2} & 0 & \check{I}_{m} \\ 0 & \check{I}_{m}\herm & 0\end{bmatrix},
\end{equation}
where \(\check{I}_{m}=[I_{m},0]\herm\in\mathbb{C}^{n\times m}\);
\(Z=[Z_{1}\herm,Z_{2}\herm]\herm\) with \(Z_{2}\in\mathbb{C}^{p\times n}\).
Note that the approximation pertains solely to the rounding error
similarly to the arguments in Section~\ref{subsec:gmres-lse}.
In exact arithmetic~\eqref{eq:gls-case} is satisfied as an equality.

If \(n>p\), we partition \(\hat{T}\) by \(\hat{T}=[\hat{T}_{1}\herm,
\hat{T}_{2}\herm]\herm\) with \(\hat{T}_{1}\in\mathbb{C}^{(n-p)\times p}\)
and \(\hat{T}_{2}\in\mathbb{C}^{p\times p}\) being upper triangular.
Then a two-sided preconditioner is presented as
\begin{equation}
\label{eq:gls-bdpre-nlarge}
\begin{split}
M_{l} & =\diag\bigl\{\alpha^{-1/2}I_{p},\alpha^{1/2}U^{-1}Q\herm,
\alpha^{-1/2}\beta^{-1}Y\herm\hat{R}\iherm\bigr\},\\
M_{r} & =\diag\bigl\{\alpha^{-1/2}I_{p},\alpha^{1/2}QU\iherm,
\alpha^{-1/2}\beta^{-1}\hat{R}^{-1}Y\bigr\},
\end{split}
\end{equation}
where \(U\in\mathbb{C}^{n\times n}\) and
\(Y\in\mathbb{C}^{m\times m}\) are defined as
\[
U=\begin{bmatrix}I_{n-p} & \hat{T}_{1} \\ 0 & \hat{T}_{2}\end{bmatrix}
\qquad\text{and}\qquad Y=\begin{bmatrix}I_{n-p} & \hat{T}_{1}(:,1:m-n+p) \\
0 & \hat{T}_{2}(1:m-n+p,1:m-n+p)\end{bmatrix}.
\]
As in~\eqref{eq:gls-case}, one can verify that
\[
M_{l}\tilde{F}M_{r}\approx\begin{bmatrix}I_{p} & \check{Z}\herm & 0 \\
\check{Z} & 0 & \check{I}_{m} \\ 0 & \check{I}_{m}\herm & 0\end{bmatrix},
\]
where \(\check{Z}=[0,Z\herm]\herm\in\mathbb{C}^{n\times p}\).
Similarly for the GLS problem, one can prove that the condition
number of the preconditioned matrix satisfies
\[
\kappa_{\infty}(M_{l}\tilde{F}M_{r})\leq C\cdot(m+n+p),
\]
if \(\mathcal{O}(\machepsf)\kappa_{2}(W)\kappa_{2}(V)\leq1\).

The main objective of the GMRES solver is to overcome the limitation on the
condition number \(\kappa_{\infty}(\tilde{F})\) in Section~\ref{sec:accuracy}.
Compared to traditional iterative refinement techniques that have
an upper limit of \(\mathcal{O}(\machepsf^{-1})\) on the condition number
for which convergence is still achieved of, it has been shown in~\cite{CH2018}
that GMRES algorithms, if implemented carefully, can relax this constraint
to \(\mathcal{O}(\machepsf^{-1}\machepso^{-1/2})\).
The limiting accuracy of the mixed precision GMRES algorithm
is no different to the classical iterative refinement ones
in Theorems~\ref{thm:lse-limit} and~\ref{thm:gls-limit}.
A specific example of that can be found in Table~\ref{tab:gmres}
by employing IEEE half/single/double precisions.

\begin{table}[!tb]
\begin{spacing}{1.3}
\caption{Comparison of results for classical and GMRES iterative refinement.}
\label{tab:gmres}
\medskip\centering
\begin{tabular}{cccccc}\hline
Method  &  \(\machepsf\) & \(\machepso\) & \(\machepsr\)
& Upper bound on \(\kappa_{\infty}(\tilde{F})\)  &  Forward error \\ \hline
classical &   half   &  single  &  double  &  \(\approx10^{4}\)   &  single \\
classical &  single  &  double  &  double
& \(\approx10^{8}\)  &  \(\cond(\tilde{F},\tilde{u})\cdot10^{-16}\) \\
classical &  single  &  double  &   quad   &  \(\approx10^{8}\)   &  double \\
  GMRES   &   half   &  single  &  double  &  \(\approx10^{8}\)   &  single \\
  GMRES   &  single  &  double  &  double  &  \(\approx10^{16}\)
& \(\cond(\tilde{F},\tilde{u})\cdot10^{-16}\) \\
  GMRES   &  single  &  double  &   quad   &  \(\approx10^{16}\)  &  double \\ 
\hline
\end{tabular}
\end{spacing}
\end{table}
\section{Numerical experiments}
\label{sec:experiments}
In this section, we present numerical results
of our mixed precision LSE and GLS algorithms.
In our tests, we set the combinations of precisions
as (\(\machepsf\), \(\machepss\), \(\machepso\), \(\machepsr\))
= (single, single, single, double) in IEEE standard.
All experiments are performed on a Linux server equipped with two sixteen-core
Intel Xeon Gold 6226R 2.90 GHz CPUs with 1024 GB of main memory.

\subsection{Experiment settings}
For the LSE problem~\eqref{problem:lse},
the following four algorithms are tested:
    
\begin{enumerate}
\item  \texttt{DGGLSE} [the reference algorithm]:
the LAPACK subroutine using the nullspace approach via
the GRQ factorization performed totally in double precision.
\item  \textbf{MPLSE}: mixed precision LSE algorithm with classical
iterative refinement presented in Algorithm~\ref{alg:mplse}.
\item  \textbf{MPLSE-GMRES-Left}: mixed precision LSE algorithm using GMRES
with the left preconditioner~\eqref{eq:lse-leftpre} as iterative refinement.
\item  \textbf{MPLSE-GMRES-BD}: mixed precision LSE algorithm using GMRES
with the block-diagonal split preconditioner~\eqref{eq:lse-bdpre-mlarge}
or~\eqref{eq:lse-bdpre-nlarge} as iterative refinement.
\end{enumerate}

Similarly, three algorithms are tested to
solve the GLS problem~\eqref{problem:gls}:

\begin{enumerate}
\item  \texttt{DGGGLM} [the reference algorithm]:
the LAPACK subroutine using Paige's algorithm via
the GQR factorization performed totally in double precision.
\item  \textbf{MPGLS}: mixed precision GLS algorithm with classical
iterative refinement presented in Algorithm~\ref{alg:mpgls}.
\item  \textbf{MPGLS-GMRES-BD}: mixed precision GLS algorithm using GMRES
with the block-diagonal split preconditioner~\eqref{eq:gls-bdpre-plarge}
or~\eqref{eq:gls-bdpre-nlarge} as iterative refinement.
\end{enumerate}

As stated in the beginning of this section, our experiments
utilize two precisions: IEEE double and single precision.
When employing a lower precision, it is often necessary to apply
appropriate scaling to prevent overflow during the conversion of
data from a higher precision to a lower precision in practice.

In our implementation, the GMRES algorithm
is conducted in double precision entirely.
Carson and Dau\v{z}ickait\.{e} discussed using lower precision
in GMRES-based iterative refinement in~\cite{CD2024}.
However, for dense matrices, there is no need for that,
as the most time-consuming part, matrix--vector multiplication,
needs to be performed in the working precision regardless.
The parameter \texttt{tol} is set to \(10^{-13}\)
for the stopping criteria of iterative refinement.

\subsection{Tests for the LSE problem}
\label{subsec:test-lse}
For the LSE problem, cases \(m>n\gg p\) and \(m\gg n>p\)
are often of interest in practice.
In case \(m>n\gg p\), we choose \(n\in\{1024,2048,3072\}\) in our experiments,
and set \((m,p)=(8n,n/32)\), \((8n,n/64)\), \((10n,n/32)\), \((10n,n/64)\),
\((12n,n/32)\), and \((12n,n/64)\) to generate different test matrices with
fixed condition number \(\kappa_{2}\big([A\herm,B\herm]\big)=10^{3}\),
\(10^{5}\), \(10^{7}\), \(10^{9}\), respectively.
In case \(m\gg n>p\), we choose \(n\in\{128,192,256\}\),
and set \((m,p)=(128n,n/8)\), \((128n,n/4)\), \((128n,n/2)\),
\((256n,n/8)\), \((256n,n/4)\), and \((256n,n/2)\) with fixed
condition number \(\kappa_{2}\big([A\herm,B\herm]\big)=10^{5}\).
The matrix \([A\herm,B\herm]\) with a specified condition number is
generated through \([A\herm,B\herm]=W_{1}\Sigma W_{2}\), where \(\Sigma\) is
a diagonal matrix with the given condition number produced by \texttt{DLATM1}
from LAPACK, and \(W_{1}\) and \(W_{2}\) are the unitary factors
from the QR factorization of randomly generated matrices.
The parameters are set as \(\alpha=\norm{r_{0}}\) and \(\beta=1\).

In each plot, we show the relative run time, i.e., the ratio of the wall clock
time of a solver over that of \texttt{DGGLSE}, and use labels `\texttt{sgglse}',
`init \(v\) and \(r\)', `residual', `correction system' to
represent the different components of classical iterative
refinement in Algorithm~\ref{alg:mplse}, respectively:
computing the initial guess by \texttt{SGGLSE},
computing the initial guess \(v\) and \(r\),
computing the residual in iterative refinement,
solving the correction system in iterative refinement.
We use labels `GMRES-Left', `GMRES-BD', respectively, to represent
iterative refinement by GMRES with the left and the block-diagonal
split preconditioner suggested in Section~\ref{sec:gmres}.
Other components of the algorithms are labelled as `others'.

\subsubsection{Case \(m>n\gg p\)}
\begin{figure}[!tb]\centering
\includegraphics[width=0.8\textwidth]{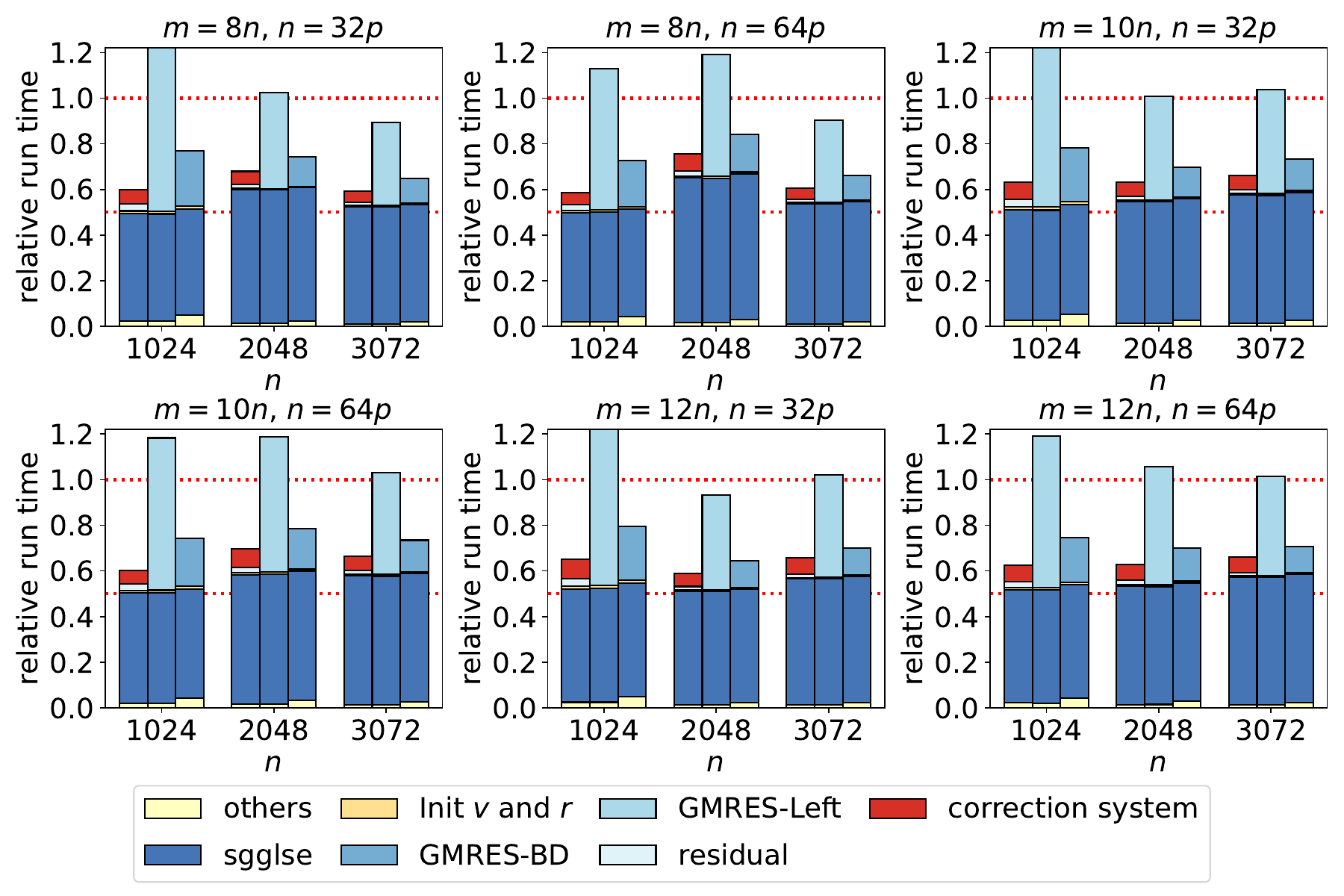}
\caption{Relative run time of MPLSE, MPLSE-GMRES-Left,
and MPLSE-GMRES-BD compared to \texttt{DGGLSE} for matrices
with \(\kappa_{2}\bigl([A\herm,B\herm]\bigr)=10^{3}\).
For each matrix, the three columns from left to right represent the results
of MPLSE, MPLSE-GMRES-Left, and MPLSE-GMRES-BD, respectively.}
\label{fig:performance-lse-3}
\end{figure}

\begin{figure}[!tb]\centering
\includegraphics[width=0.8\textwidth]{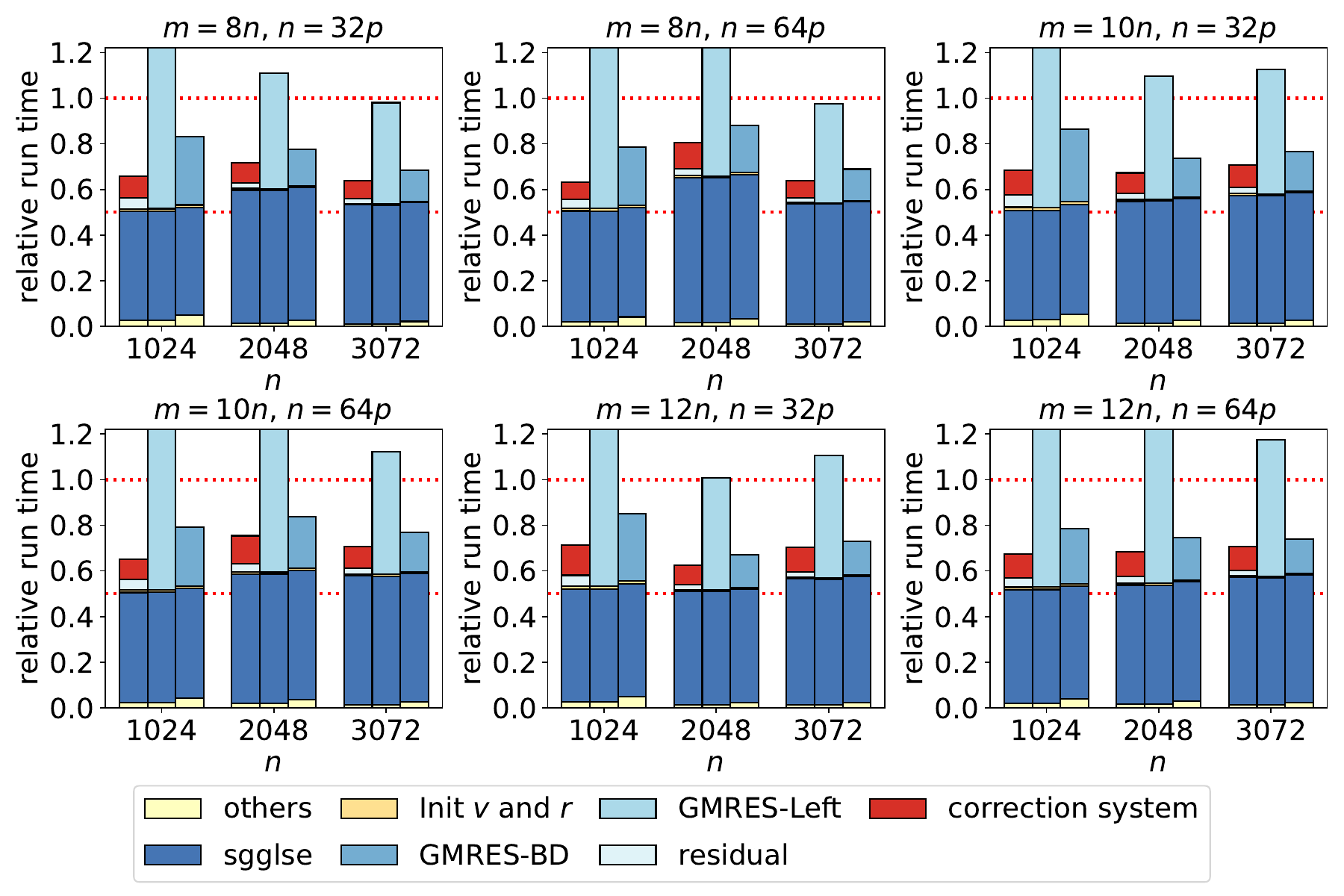}
\caption{Relative run time of MPLSE, MPLSE-GMRES-Left,
and MPLSE-GMRES-BD compared to \texttt{DGGLSE} for matrices
with \(\kappa_{2}\bigl([A\herm,B\herm]\bigr)=10^{5}\).
For each matrix, the three columns from left to right represent the results
of MPLSE, MPLSE-GMRES-Left, and MPLSE-GMRES-BD, respectively.}
\label{fig:performance-lse-5}
\end{figure}

\begin{figure}[!tb]\centering
\includegraphics[width=0.8\textwidth]{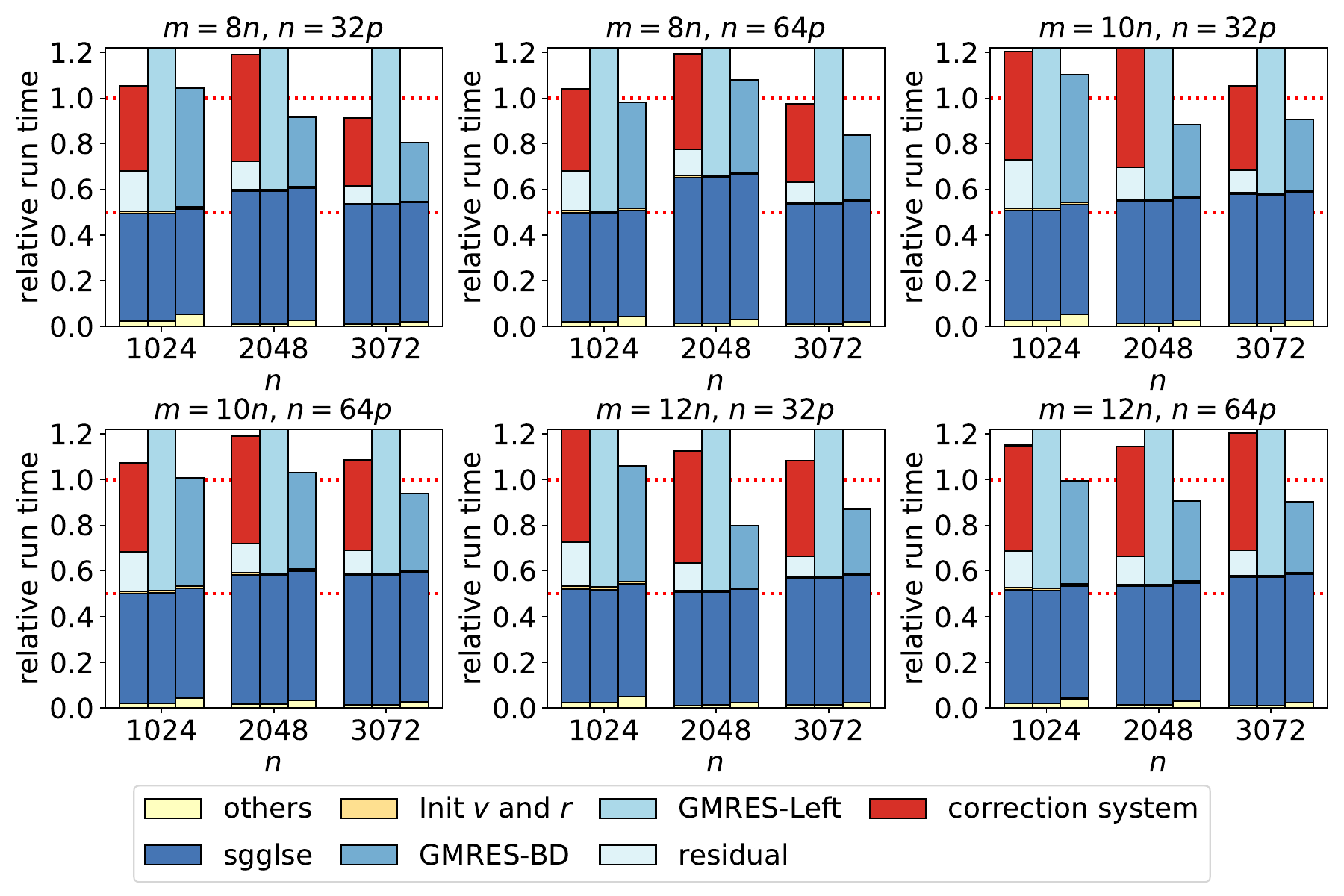}
\caption{Relative run time of MPLSE, MPLSE-GMRES-Left,
and MPLSE-GMRES-BD compared to \texttt{DGGLSE} for matrices
with \(\kappa_{2}\bigl([A\herm,B\herm]\bigr)=10^{7}\).
For each matrix, the three columns from left to right represent the results
of MPLSE, MPLSE-GMRES-Left, and MPLSE-GMRES-BD, respectively.}
\label{fig:performance-lse-7}
\end{figure}

We first focus on case \(m>n\gg p\).
Table~\ref{tab:accuracy-lse} presents the accuracy
of the results and the iteration counts.
Figures~\ref{fig:performance-lse-3}--\ref{fig:performance-lse-7} illustrate
the relative run time for matrices with various condition numbers.
It is important to note that Table~\ref{tab:accuracy-lse} displays
the results for a single test matrix of consistent size across varying
condition numbers, as the accuracy and iteration counts are predominantly
affected by the condition number rather than the matrix dimension.

\begin{table}[!tb]
\begin{spacing}{1.3}
\caption{Accuracy and iteration counts of MPLSE,
MPLSE-GMRES-Left and MPLSE-GMRES-BD.
`err-1', `err-2', and `iter' denote, respectively,
\(\norm{Bx-d}/\big(\normfro{B}\norm{x}+\norm{d}\big)\), \(\bigabs{\norm{Ax-b}
/\norm{Ax_{\mathtt{dgglse}}-b}-1}\), which is related to \(\kappa_{2}
\bigl([A\herm,B\herm]\bigr)\), and the number of iterations.}
\label{tab:accuracy-lse}
\medskip\centering
\begin{tabular}{c|c|cccc}\hline
\multicolumn{2}{c|}{\(\kappa_{2}\big([A\herm,B\herm]\big)\)} & \(\kappa=10^{3}\)
& \(\kappa= 10^{5}\)  &  \(\kappa=10^{7}\) & \(\kappa=10^{9}\) \\ \hline
\multirow{3}{*}{MPLSE}
&  err-1  &  \(3.3\cdot 10^{-17}\)  &  \(2.0\cdot 10^{-16}\)
&  \(2.2\cdot 10^{-14}\) &  \(2.6\cdot 10^{-10}\) \\
&  err-2  &  \(2.9\cdot 10^{-16}\)  &  \(5.8\cdot 10^{-14}\)
&  \(9.9\cdot 10^{-11}\) &  \(7.6\cdot 10^{-2}\) \\
&   iter  &  \(2\)  &  \(3\)  &  \(12\) & diverged \\ \hline
\multirow{3}{*}{MPLSE-GMRES-Left}
&  err-1  &  \(2.7\cdot 10^{-17}\)  &  \(5.0\cdot 10^{-16}\)
&  \(2.3\cdot 10^{-17}\) &  \(7.2\cdot 10^{-16}\) \\
&  err-2  &  \(4.3\cdot 10^{-16}\)  &  \(1.6\cdot 10^{-13}\)
&  \(4.1\cdot 10^{-12}\) &  \(2.3\cdot 10^{-9}\) \\
&   iter  &  \(4\)  &  \(5\)  &  \(15\) & \(145\) \\ \hline
\multirow{3}{*}{MPLSE-GMRES-BD}
&  err-1  &  \(3.4\cdot 10^{-17}\)  &  \(4.0\cdot 10^{-16}\)
&  \(1.1\cdot 10^{-14}\) &  \(3.7\cdot 10^{-17}\) \\
&  err-2  &  \(7.2\cdot 10^{-16}\)  &  \(1.7\cdot 10^{-13}\)
&  \(5.6\cdot 10^{-11}\) &  \(3.9\cdot 10^{-10}\) \\
&   iter  &  \(11\)  &  \(14\)  &  \(30\) & \(243\) \\ \hline
\end{tabular}
\end{spacing}
\end{table}

From Table~\ref{tab:accuracy-lse} and
Figures~\ref{fig:performance-lse-3}--\ref{fig:performance-lse-7},
we observe that in comparison with the left preconditioner,
the GMRES-based algorithm with the block-diagonal split
preconditioner requires more iterations to converge.
However, the left preconditioner costs more execution time due to the higher
computational cost of applying~\eqref{eq:lse-leftpre} to a vector.
This indicates that it is preferable to use the block-diagonal split
preconditioner~\eqref{eq:lse-bdpre-mlarge} or~\eqref{eq:lse-bdpre-nlarge}
for GMRES-based iterative refinement.

For test matrices with relatively small condition numbers, as shown in
Figures~\ref{fig:performance-lse-3} and~\ref{fig:performance-lse-5},
MPLSE is the most efficient among the three mixed precision LSE algorithms,
which saves around \(30\%\)--\(40\%\) of the execution time
over \texttt{DGGLSE}. This is roughly a speedup around 1.5.

As the condition number \(\kappa_{2}\bigl([A\herm,B\herm]\bigr)\)
becomes \(10^{7}\), which is close to \(\machepsf^{-1}\),
all three mixed precision LSE algorithms require more iterations
and longer time to converge; see Figure~\ref{fig:performance-lse-7}.
This is due to the reduced accuracy of the initial guess and GRQ factorization,
which also affects the quality of the preconditioner.
Our analysis indicates that MPLSE is not suitable for scenarios
where condition numbers exceed approximately \(\machepsf^{-1}\).
As demonstrated by the experiments in Table~\ref{tab:accuracy-lse}
with \(\kappa=10^{9}\), MPLSE fails to converge in this case.
In Figure~\ref{fig:performance-lse-7}, we see that
MPLSE-GMRES-BD remains approximately \(10\%\)--\(20\%\) faster
than \texttt{DGGLSE}, a speedup between 1.11 and 1.25.
However, that will diminish rapidly as the condition number increases further.

\subsubsection{Case \(m\gg n>p\)}
In case \(m\gg n>p\), comparison among MPLSE, MPLSE-GMRES-Left,
and MPLSE-GMRES-BD is similar to the case \(m>n\gg p\),
which we will no longer reiterate.

\begin{figure}[!tb]\centering
\includegraphics[width=0.8\textwidth]{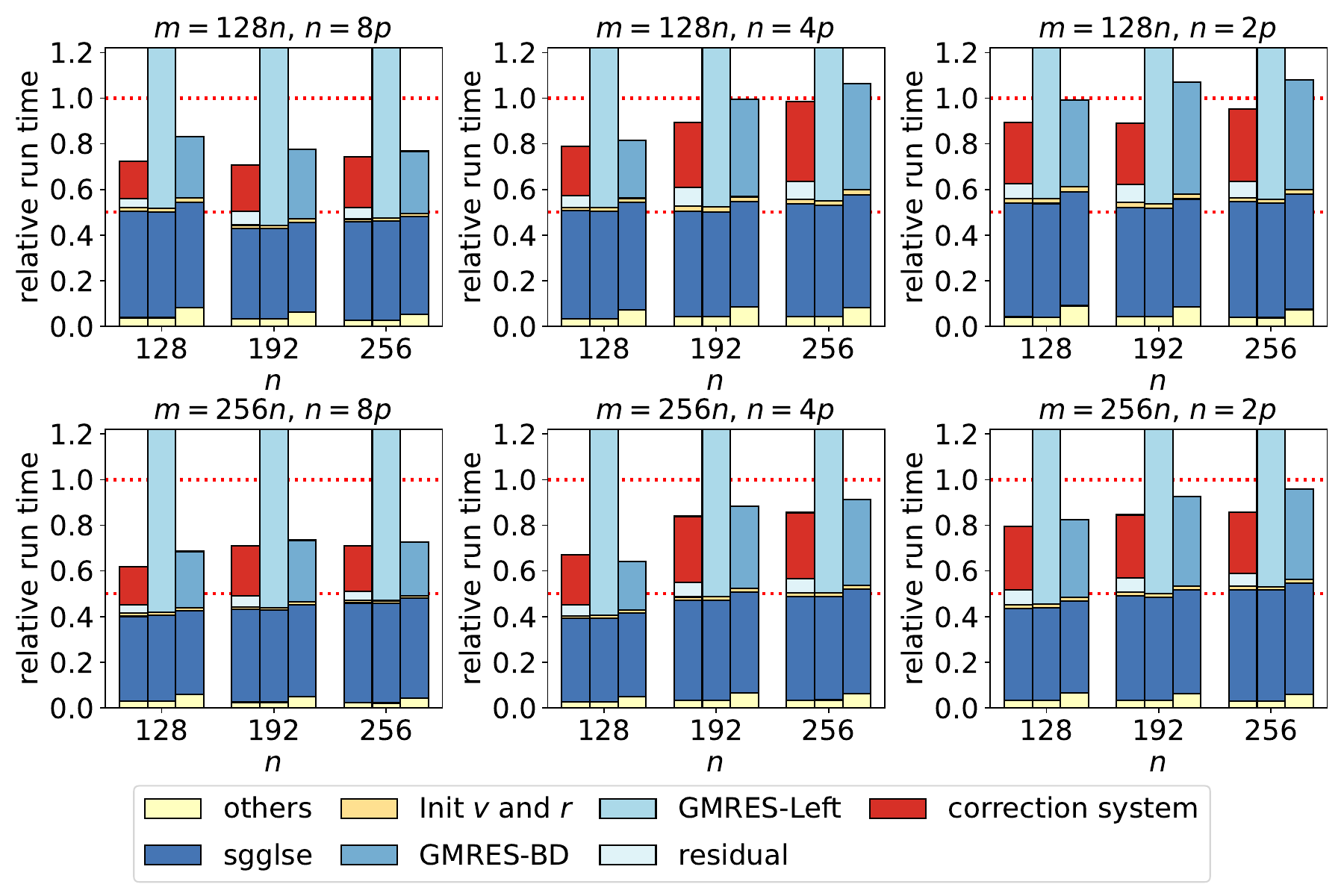}
\caption{Relative run time of MPLSE, MPLSE-GMRES-Left,
and MPLSE-GMRES-BD compared to \texttt{DGGLSE} for matrices
with \(\kappa_{2}\bigl([A\herm,B\herm]\bigr)=10^{5}\).
For each matrix, the three columns from left to right represent the
results of MPLSE, MPLSE-GMRES-Left, and MPLSE-GMRES-BD, respectively.}
\label{fig:performance-lse-5-case2}
\end{figure}

However, the overall speedup in this scenario is not as significant
as in the case \(m>n\gg p\), particularly when \(n\) is close to \(p\).
As depicted in Figure~\ref{fig:performance-lse-5-case2}, employing iterative
refinement-based mixed precision techniques may reduce computational time
by approximately \(40\%\) when \(n=8p\), a speed up around 1.66.
However, when \(n=4p\) or \(n=2p\),
the time savings become progressively smaller.
This can be explained from the complexity perspective
discussed in Section~\ref{subsec:mpalgo-lse}.
When \(n\) is relatively small, \(\mathcal{O}(mn^{2})\) becomes
close to \(\mathcal{O}(mn)\), resulting in a larger proportion of
execution time being consumed by iterative refinement
and thus diminishing the performance of the algorithm.

\subsection{Tests for the GLS problem}
\label{subsec:test-gls}
For the GLS problem, we choose the dimension \(n\in\{1024,2048,3072\}\) in
our experiments, and set \((p,m)=(8n,n/32)\), \((8n,n/64)\), \((10n,n/32)\),
\((10n,n/64)\), \((12n,n/32)\), and \((12n,n/64)\) to generate different test matrices with fixed condition number \(\kappa_{2}\bigl([W,V]\bigr)=10^{3}\),
\(10^{5}\), \(10^{7}\), \(10^{9}\), respectively.
Note that we do not test GMRES with the left
preconditioner~\eqref{eq:gls-leftpre} for the GLS problem,
since it has already been shown in Section~\ref{subsec:test-lse} that
it is preferable to use the block-diagonal split preconditioner.
Similarly, the parameters are set as \(\alpha=\norm{y_{0}}\) and \(\beta=1\).

Similarly to Section~\ref{subsec:test-lse}, in each plot we show
the relative run time and use labels `\texttt{sggglm}', `init \(z\)',
`residual', `correction system' to represent the different components
of classical iterative refinement in Algorithm~\ref{alg:mpgls}, respectively.
We use labels `GMRES-BD' to represent iterative refinement by GMRES with
the block-diagonal split preconditioner suggested in Section~\ref{sec:gmres}.
Other components of the algorithms are labelled as `others'.

\begin{figure}[!tb]\centering
\includegraphics[width=0.8\textwidth]{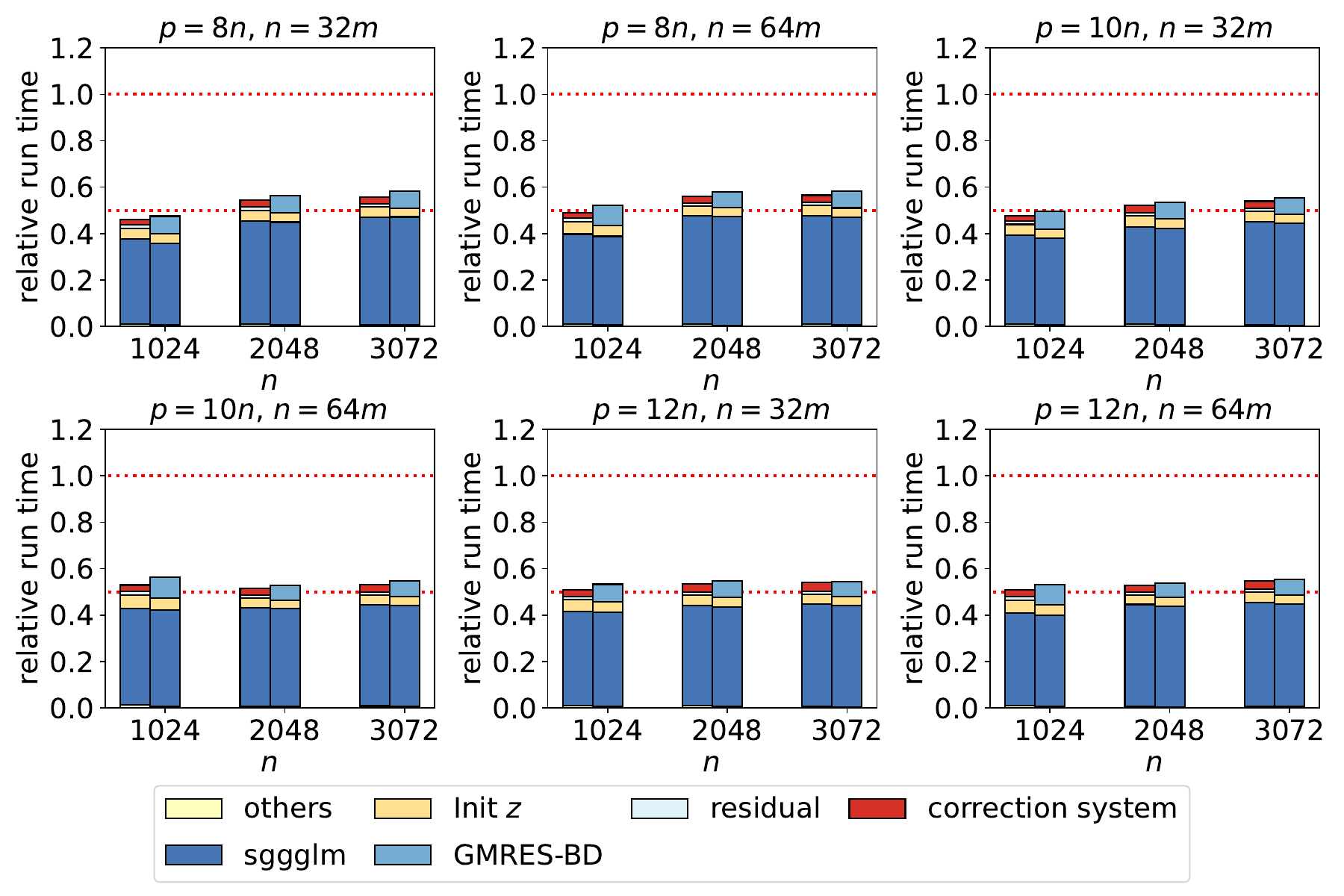}
\caption{Relative run time of MPGLS and MPGLS-GMRES-BD compared to
\texttt{DGGGLM} for matrices with \(\kappa_{2}\bigl([W,V]\bigr)=10^{3}\).
For each matrix, the two columns from left to right represent
the result of MPGLS and MPGLS-GMRES-BD, respectively.}
\label{fig:performance-gls-3}
\end{figure}

\begin{figure}[!tb]\centering
\includegraphics[width=0.8\textwidth]{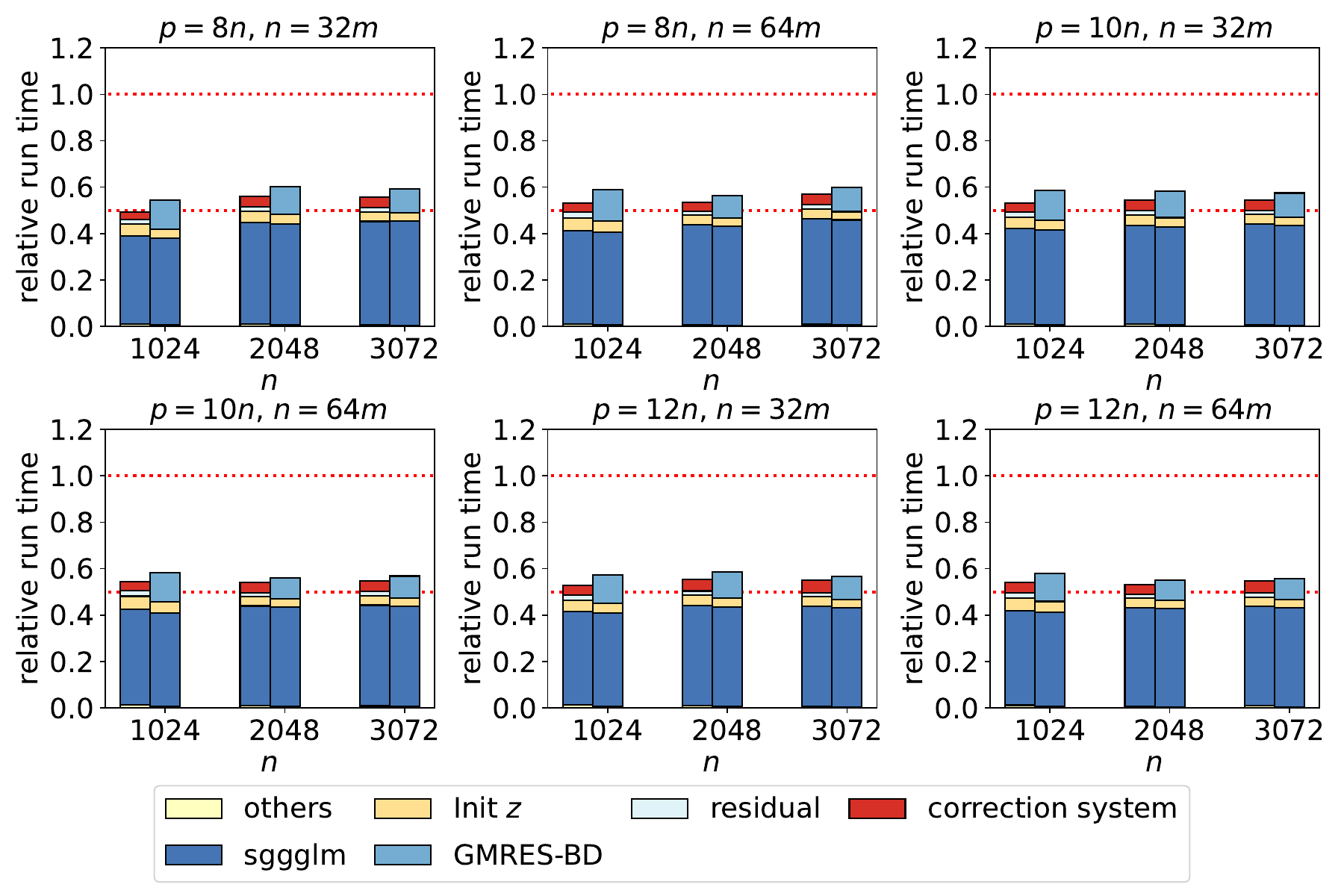}
\caption{Relative run time of MPGLS and MPGLS-GMRES-BD compared to
\texttt{DGGGLM} for matrices with \(\kappa_{2}\bigl([W,V]\bigr)=10^{5}\).
For each matrix, the two columns from left to right represent
the result of MPGLS and MPGLS-GMRES-BD, respectively.}
\label{fig:performance-gls-5}
\end{figure}

\begin{figure}[!tb]\centering
\includegraphics[width=0.8\textwidth]{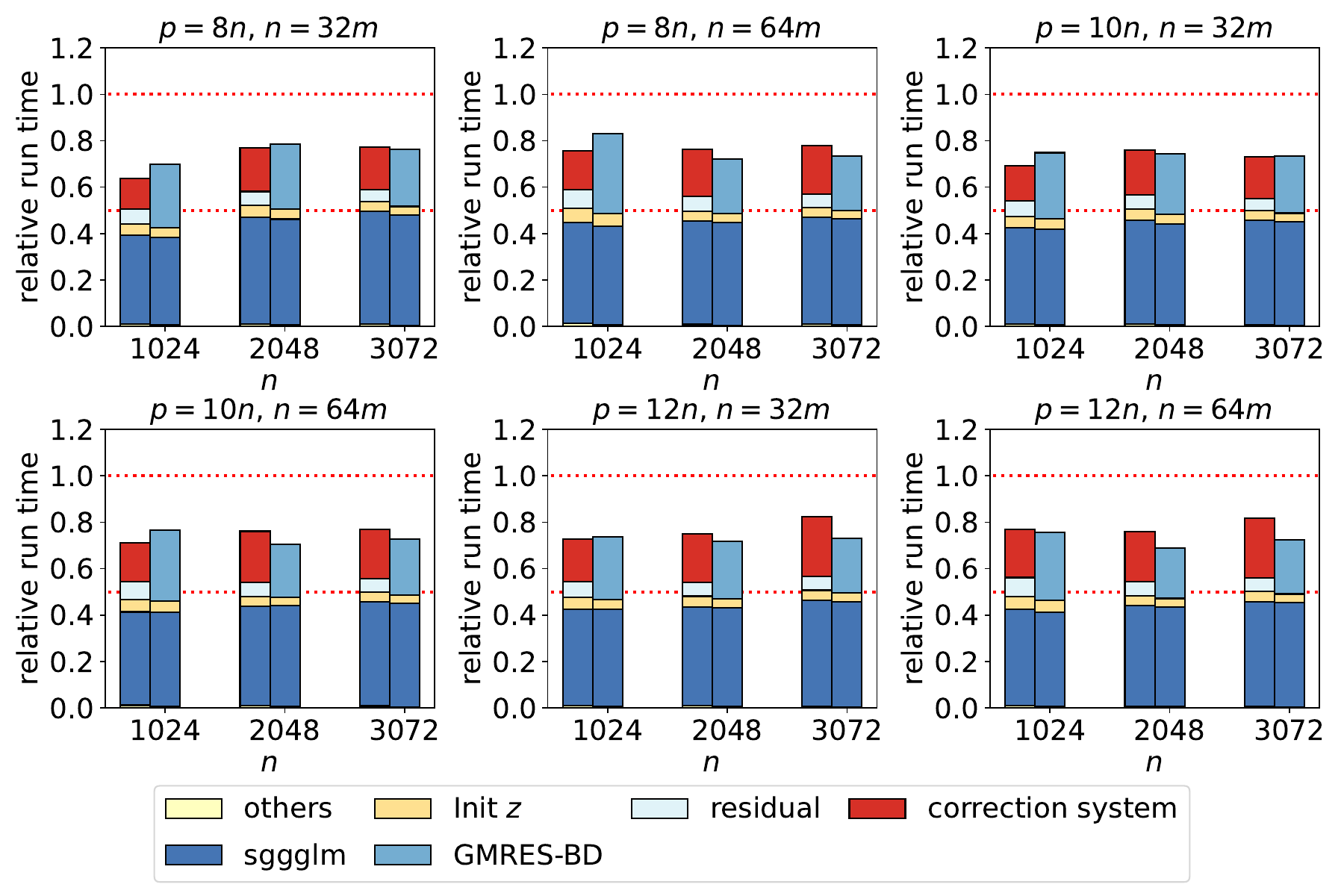}
\caption{Relative run time of MPGLS and MPGLS-GMRES-BD compared to
\texttt{DGGGLM} for matrices with \(\kappa_{2}\bigl([W,V]\bigr)=10^{7}\).
For each matrix, the two columns from left to right represent
the result of MPGLS and MPGLS-GMRES-BD, respectively.}
\label{fig:performance-gls-7}
\end{figure}

Table~\ref{tab:accuracy-gls} presents the accuracy of the results
and the iteration counts for a single test matrix of consistent
size across varying condition numbers.
In Figures~\ref{fig:performance-gls-3}--\ref{fig:performance-gls-7},
we present the relative run time for matrices with fixed condition numbers
\(\kappa_{2}\bigl([W,V]\bigr)=10^{3}\), \(10^{5}\), \(10^{7}\).
For test matrices with relatively small condition numbers, both algorithms
reduce the computational time by a factor around \(40\%\)--\(50\%\), i.e.,
a speedup between \(1.66\) and \(2.00\), with MPGLS being slightly faster.
The analysis shows that MPGLS fails when condition numbers exceed
approximately~\(\machepsf^{-1}\), which is also confirmed by the results
presented in Table~\ref{tab:accuracy-gls} with \(\kappa=10^{9}\).
As demonstrated in Table~\ref{tab:accuracy-gls}, MPGLS-GMRES-BD can still solve
these relatively ill-conditioned problems, likely with no performance gain.

\begin{table}[!htb]
\begin{spacing}{1.3}
\caption{Accuracy and iteration counts of MPGLS and MPGLS-GMRES-BD.
`err-1', `err-2', and `iter' denote, respectively, \(\norm{Wx+Vy-d}
/\bigl(\normfro{W}\norm{x}+\normfro{V}\norm{y}+\norm{d}\bigr)\),
\(\bigabs{\norm{y}/\norm{y_{\mathtt{dggglm}}}-1}\), which is related to
\(\kappa_{2}\bigl([W,V]\bigr)\), and the number of iterations.}
\label{tab:accuracy-gls}
\medskip\centering
\begin{tabular}{c|c|cccc}\hline
\multicolumn{2}{c|}{\(\kappa_{2}\big([W,V]\big)\)}
&  \(\kappa=10^{3}\)  &  \(\kappa=10^{5}\)
&  \(\kappa=10^{7}\)  &  \(\kappa=10^{9}\) \\ \hline
\multirow{3}{*}{MPGLS}
&  err-1  &  \(2.0\cdot 10^{-17}\)  &  \(5.0\cdot 10^{-16}\)
&  \(9.4\cdot 10^{-15}\) & \(8.7\cdot 10^{-10}\) \\
&  err-2  &  \(4.1\cdot 10^{-15}\)  &  \(1.0\cdot 10^{-11}\)
&  \(7.2\cdot 10^{-8}\)  & \(6.8\cdot 10^{-1}\) \\
&   iter  &  \(3\)  &  \(4\)  &  \(13\) & diverged \\ \hline
\multirow{3}{*}{MPGLS-GMRES-BD}
&  err-1  &  \(8.9\cdot 10^{-15}\)  &  \(2.9\cdot 10^{-14}\)
&  \(4.1\cdot 10^{-13}\) & \(6.7\cdot 10^{-11}\) \\
&  err-2  &  \(2.3\cdot 10^{-15}\)  &  \(2.4\cdot 10^{-13}\)
&  \(8.7\cdot 10^{-12}\) & \(1.1\cdot 10^{-8}\) \\
&   iter  &  \(11\)  &  \(18\)  &  \(47\) & \(297\) \\ \hline
\end{tabular}
\end{spacing}
\end{table}

\section{Conclusions}
\label{sec:conclusions}
In this paper, we propose iterative refinement-based mixed precision algorithms
for two specific variants of least squares problems---least squares problems
with linear equality constraints and generalized least squares problems.

For each problem, we prove that mixed precision classical
iterative refinement algorithms converge for reasonably well-conditioned
matrices that satisfy \(\machepsf\kappa_{\infty}(\tilde{F})\lesssim 1\),
with \(\tilde{F}\) being the augmented matrix.
To overcome the limitation on the condition number, we further develop
a GMRES-based iterative refinement solver along with two preconditioners%
---a left preconditioner and a block-diagonal split preconditioner.
Both classical and GMRES-based mixed precision iterative
refinement algorithms are capable of refining the solution to
the working precision level in a backward stable manner.
We also show that employing a higher precision in the computation of residuals
allows mixed precision algorithms to further improve the accuracy of the
computed results to the working precision level in a forward stable manner.

Numerical experiments indicate that when the condition number
\(\kappa_{\infty}(\tilde{F})\) is not too large, our mixed precision
algorithms enhance the performance significantly.
Classical iterative refinement proves to be the fastest of the tested
algorithms and only costs around \(60\%\) of the time of
the LAPACK subroutines \texttt{DGGLSE} and \texttt{DGGGLM}.
This shows that the mixed precision classical iterative
refinement algorithm is preferably suitable in numerical linear
algebra libraries for reasonably well-conditioned matrices.
As the condition number increases, the classical iterative
refinement algorithm may potentially fail to converge,
and the mixed precision GMRES-based iterative refinement algorithm with
the block-diagonal split preconditioner emerges as a favourable option.
None of these iterative refinement based mixed precision algorithms
achieves any acceleration for problems that are overly ill-conditioned.
In particular, in our experiments we only observe speedup
for linear systems with condition number up to \(10^{7}\),
i.e., roughly \(\mathcal{O}(\machepsf^{-1})\).
\bibliographystyle{plainurl}
\bibliography{reference}
\appendix
\section{Proof of Theorem~\ref{thm:lse}}
\label{sec:appendix-proof}
\begin{proof}
By Lemma~\ref{lemma:grq}, our goal shifts to prove that there exist
\(\Delta\breve{A}_{1}\), \(\Delta\breve{A}_{2}\),
\(\Delta\breve{B}_{1}\), \(\Delta\breve{B}_{2}\),
and \(\Delta t_{1}\), \(\Delta t_{2}\), \(\Delta t_{3}\) such that
\begin{equation}
\label{eq:modified-equation-lse}
\begin{bmatrix}
I & 0 & Z\hat{T}Q+\Delta\breve{A}_{1} \\ 0 & 0 & [0,\hat{R}]Q
+\Delta\breve{B}_{1} \\ \big(Z\hat{T}Q+\Delta\breve{A}_{2}\big)\herm &
\big([0,\hat{R}]Q+\Delta\breve{B}_{2}\big)\herm & 0\end{bmatrix}\begin{bmatrix}
\fl(\Delta r) \\ -\fl(\Delta v) \\ \fl(\Delta x)\end{bmatrix}=\begin{bmatrix}
f_{1}+\Delta t_{1} \\ f_{2}+\Delta t_{2} \\ f_{3}+\Delta t_{3}\end{bmatrix}.
\end{equation}
Then we conclude the proof by defining \(\Delta B_{i}=\Delta\breve{B}_{i}
+\Delta E_{1}\) and \(\Delta A_{i}=\Delta\breve{A}_{i}+\Delta E_{2}\).
As discussed in Section~\ref{sec:algorithm},
the solution to~\eqref{eq:modified-equation-lse}
is obtained by addressing~\eqref{eq:sum-lse}.
Consequently, we shall analyse the rounding errors encountered
in solving~\eqref{eq:sum-lse} sequentially following the
computational steps outlined in Algorithm~\ref{alg:lse-solver}.

In Algorithm~\ref{alg:lse-solver}, we first
compute \(u=Qf_{3}\) and \(w=Z\herm f_{1}\).
By Assumption~\ref{assump:unmqr}, there exist \(\Delta s_{1}\in\mathbb{C}^{m}\)
and \(\Delta s_{3}\in\mathbb{C}^{n}\) such that
\[
\hat{w}=Z\herm(f_{1}+\Delta s_{1}),\qquad\hat{u}=Q(f_{3}+\Delta s_{3}),
\]
where \(\norm{\Delta s_{1}}\leq\epsunmqr^{s}(m)\norm{f_{1}}\)
and \(\norm{\Delta s_{3}}\leq\epsunmqr^{s}(n)\norm{f_{3}}\).

\medskip
We next consider the rounding errors in solving~\eqref{eq:y2} and~\eqref{eq:q1}.
For the triangular systems \(Ry_{2}=f_{2}\) and \(T_{11}\herm q_{1}=u_{1}\), by
Assumption~\ref{assump:trsv}, there exist \(\Delta E_{3}\) and \(\Delta E_{4}\)
such that the computed quantities \(\hat{y}_{2}\) and \(\hat{q}_{1}\) satisfy
\begin{equation}
\label{eq:error-y2q1}
(\hat{R}+\Delta E_{3})\hat{y}_{2}=f_{2},\qquad
(\hat{T}_{11}+\Delta E_{4})\herm\hat{q}_{1}=\hat{u}_{1},
\end{equation}
where \(\normfro{\Delta E_{3}}\leq\epstrsv^{s}(p)\normfro{\hat{R}}\) and
\(\normfro{\Delta E_{4}}\leq\epstrsv^{s}(n-p)\normfro{\hat{T}_{11}}\).

\medskip
Then we analyse~\eqref{eq:y1} and~\eqref{eq:q2}.
By Assumption~\ref{assump:gemv}, there exist
\(\Delta E_{5}\) and \(\Delta E_{6}\) such that
\begin{align*}
\fl(T_{12}y_{2})=(\hat{T}_{12}+\Delta E_{5})\hat{y}_{2},\qquad
\fl(T_{22}y_{2})=(\hat{T}_{22}+\Delta E_{6})\hat{y}_{2},
\end{align*}
where \(\normfro{\Delta E_{5}}\leq\epsgemv^{s}(p)\normfro{\hat{T}_{12}}\)
and \(\normfro{\Delta E_{6}}\leq\epsgemv^{s}(p)\normfro{\hat{T}_{22}}\).
Furthermore, from the standard rounding error analysis we derive that there
exist \(\epsilon_{1}\), \(\epsilon_{2}\), and matrices \(J_{1}\), \(J_{2}\)
such that the right-hand sides of~\eqref{eq:y1} and~\eqref{eq:q2} satisfy
\begin{align}
\hat{q}_{2} & =(I+\epsilon_{1}J_{1})\big(\hat{w}_{2}
-\fl(T_{22}y_{2})\big),\label{eq:error-q2} \\
\fl(w_{1}-q_{1}-T_{12}y_{2}) & =(I+\epsilon_{2}J_{2})\big(\hat{w}_{1}
-\hat{q}_{1}-\fl(T_{12}y_{2})\big),\label{eq:e7-1}
\end{align}
with \(\abs{\epsilon_{1}}\leq\machepss\), \(\abs{\epsilon_{2}}
\leq2\machepss+\machepss^{2}\), and \(\norm{J_{i}}=1\).
Then by Lemma~\ref{lemma:error-tri},
we summarize the errors from the right-hand side vector summations
and the triangular solver in \(\Delta E_{7}\), i.e.,
\begin{equation}
\label{eq:e7-2}
(\hat{T}_{11}+\Delta E_{7})\hat{y}_{1}=\hat{w}_{1}-\hat{q}_{1}
-\fl(T_{12}y_{2})=(I+\epsilon_{2}J_{2})^{-1}\fl(w_{1}-q_{1}-T_{12}y_{2}),
\end{equation}
where
\[
\normfro{\Delta E_{7}}\leq\frac{2\machepss+\machepss^{2}
+\epstrsv^{s}(n-p)}{1-2\machepss-\machepss^{2}}\normfro{\hat{T}_{11}}
\leq 4\epstrsv^{s}(n-p)\normfro{\hat{T}_{11}}.
\]

Substitute~\eqref{eq:error-y2q1},~\eqref{eq:error-q2},
and~\eqref{eq:e7-2} into \(\Delta r=Zq\) and \(\Delta x=Q\herm y\).
From Lemma~\ref{lemma:error-unitary} and \(\machepsf\ll1\),
we derive that there exists \(\delta r\) such that
\[
\fl(\Delta r)=Z\begin{bmatrix}\hat{q}_{1} \\
\hat{w}_{2}-\fl(T_{22}y_{2})\end{bmatrix}+\delta r,
\qquad\norm{\delta r}\leq3\epsunmqr^{s}(m)\norm{\fl(\Delta r)}.
\]
For \(\Delta x=Q\herm y\), from Assumption~\ref{assump:unmqr}
we know that there exists \(\delta x\) such that
\[
\fl(\Delta x)=Q\herm\hat{y}+\delta x,\qquad
\norm{\delta x}\leq\epsunmqr^{s}(n)\norm{\fl(\Delta x)}.
\]

Finally, for the triangular system~\eqref{eq:deltav}, i.e.,
\(R\herm\Delta v=T_{12}\herm q_{1}+T_{22}\herm q_{2}-u_{2}\),
we treat it in two steps similarly to the analysis of~\eqref{eq:y1}.
First there exist \(\Delta E_{8}\) and \(\Delta E_{9}\) such that
\[
\fl(T_{12}\herm q_{1})=(T_{12}+\Delta E_{8})\herm q_{1},\qquad
\fl(T_{22}\herm q_{2})=(T_{22}+\Delta E_{9})\herm q_{2},
\]
where \(\normfro{\Delta E_{8}}\leq\epsgemv^{s}(n-p)\normfro{\hat{T}_{12}}\)
and \(\normfro{\Delta E_{9}}\leq\epsgemv^{s}(m-n+p)\normfro{\hat{T}_{22}}\).
Then for the triangular system, as in~\eqref{eq:e7-1} and~\eqref{eq:e7-2},
there exists \(\Delta E_{10}\) such that
\[
(\hat{R}+\Delta E_{10})\herm\fl(\Delta v)
=\fl(T_{12}\herm q_{1})+\fl(T_{22}\herm q_{2})-\hat{u}_{2},
\]
where
\[
\normfro{\Delta E_{10}}\leq\frac{3\machepss+3\machepss^{2}+\machepss^{3}
+\epstrsv^{s}(p)}{1-3\machepss-3\machepss^{2}-\machepss^{3}}
\normfro{\hat{R}}\leq5\epstrsv^{s}(p)\normfro{\hat{R}}.
\]

\medskip
We summarize the backward errors in \(\Delta A_{i}\),
\(\Delta B_{i}\), and \(\Delta t_{j}\).
Let
\begin{align*}
\Delta\breve{A}_{1}=Z\begin{bmatrix}\Delta E_{7} & \Delta E_{5} \\
0 & \Delta E_{6}\end{bmatrix}Q, &\qquad
\Delta\breve{A}_{2}=Z\begin{bmatrix}\Delta E_{4} & \Delta E_{8} \\
0 & \Delta E_{9}\end{bmatrix}Q, \\
\Delta\breve{B}_{1}=[0,\Delta E_{3}]Q, &\qquad
\Delta\breve{B}_{2}=[0,\Delta E_{10}]Q, \\
\Delta t_{1} =\Delta s_{1}+\delta r+(A+\Delta A_{1})\delta x,\quad
\Delta t_{2} =(B&+\Delta B_{1})\delta x,\quad
\Delta t_{3} =\Delta s_{3}+(A+\Delta A_{2})\herm\delta r.
\end{align*}
One can verify that~\eqref{eq:modified-equation-lse} holds.
Recall from Lemma~\ref{lemma:grq} that \(\hat{T}\) and \(\hat{R}\) satisfy
\begin{align*}
\normfro{\hat{T}} & =\bignormfro{Z\herm(A+\Delta E_{2})Q\herm}
=\normfro{A+\Delta E_{2}}
\leq\normfro{A}+\normfro{\Delta E_{2}}\leq(1+\eta_{0})\normfro{A},\\
\normfro{\hat{R}} & =\bignormfro{(B+\Delta E_{1})Q\herm}
=\normfro{B+\Delta E_{1}}
\leq\normfro{B}+\normfro{\Delta E_{1}}\leq(1+\epsqr^{f}(n,p))\normfro{B}.
\end{align*}
Therefore the bounds of the matrices can be presented as
\begin{align*}
\normfro{\Delta\breve{A}_{1}} & =\Bignormfro{Z\begin{bmatrix}\Delta E_{7} &
\Delta E_{5} \\ 0 & \Delta E_{6}\end{bmatrix}Q}=\bigl(\normfro{\Delta E_{5}}^{2}
+\normfro{\Delta E_{6}}^{2}+\normfro{\Delta E_{7}}^{2}\bigr)
^{\frac{1}{2}}\leq\eta_{1}(1+\eta_{0})\normfro{A},\\
\normfro{\Delta\breve{A}_{2}} & =\Bignormfro{Z\begin{bmatrix}\Delta E_{4} &
\Delta E_{8}\\ 0 & \Delta E_{9}\end{bmatrix}Q}=\bigl(\normfro{\Delta E_{4}}^{2}
+\normfro{\Delta E_{8}}^{2}+\normfro{\Delta E_{9}}^{2}\bigr)
^{\frac{1}{2}}\leq\eta_{2}(1+\eta_{0})\normfro{A},\\
\normfro{\Delta\breve{B}_{1}} & =\bignormfro{[0,\Delta E_{3}]Q}
=\normfro{\Delta E_{3}}\leq\epstrsv^{s}(p)(1+\epsqr^{f}(n,p))\normfro{B},\\
\normfro{\Delta\breve{B}_{2}} & =\bignormfro{[0,\Delta E_{10}]Q}
=\normfro{\Delta E_{10}}\leq5\epstrsv^{s}(p)(1+\epsqr^{f}(n,p))\normfro{B},
\end{align*}
resulting in
\begin{align*}
\normfro{\Delta A_{i}}\leq\normfro{\Delta\breve{A}_{i}} & +\normfro{
\Delta E_{2}}\leq(\eta_{0}+\eta_{i}+\eta_{0}\eta_{i})\normfro{A},\\
\normfro{\Delta B_{i}}\leq\normfro{\Delta\breve{B}_{i}}
+\normfro{\Delta E_{1}} & \leq\bigl(\epsqr^{f}(n,p)+5\epstrsv^{s}(p)
+5\epsqr^{f}(n,p)\epstrsv^{s}(p)\bigr)\normfro{B},
\end{align*}
for \(i=1\), \(2\).
The bounds of the vectors are subsequently summarized as
\begin{align*}
\norm{\Delta t_{1}} & \leq\epsunmqr^{s}(m)\norm{f_{1}}
+3\epsunmqr^{s}(m)\norm{\fl(\Delta r)}+(1+\eta_{0})(1+\eta_{1})\,
\epsunmqr^{s}(n)\normfro{A}\norm{\fl(\Delta x)},\\
\norm{\Delta t_{2}} & \leq(1+\epsqr^{f}(n,p))(1+\epstrsv^{s}(p))\,
\epsunmqr^{s}(n)\normfro{B}\norm{\fl(\Delta x)},\\
\norm{\Delta t_{3}} & \leq\epsunmqr^{s}(n)\norm{f_{3}}+3(1+\eta_{0})
(1+\eta_{2})\,\epsunmqr^{s}(m)\normfro{A}\norm{\fl(\Delta r)}.\qedhere
\end{align*}
\end{proof}

\section{Rounding error regarding a `small' perturbation}
\label{sec:appendix-error}
\begin{lemma}
\label{lemma:error-tri}
Let \(L\in\mathbb{C}^{n\times n}\), \(v\in\mathbb{C}^{n}\).
For the triangular system \(Lx=v\) with \(v\) containing a small perturbation,
i.e., \(Lx=(I+\epsilon J)v\), with \(\abs{\epsilon}\leq\beta<1\),
\(\norm{J}=1\), there exists a matrix \(\Delta L\) such that
\[
(L+\Delta L)\hat{x}=v,
\]
where
\[
\normfro{\Delta L}\leq\frac{\beta+\epstrsv(n)}{1-\beta}\normfro{L}.
\]
\end{lemma}

\begin{proof}
From Table~\ref{tab:notation} we see that there
exists a matrix \(\Delta L_{1}\) such that
\[
(L+\Delta L_{1})\hat{x}=(I+\epsilon J)v,
\]
with \(\normfro{\Delta L_{1}}\leq\epstrsv(n)\cdot\normfro{L}\).
Let \(\Delta L_{2}=-\epsilon J(I+\epsilon J)^{-1}(L+\Delta L_{1})\).
Then
\begin{align*}
(L+\Delta L_{1}+\Delta L_{2})\hat{x} & =(I-\epsilon J(I+\epsilon J)^{-1})
(L+\Delta L_{1})\hat{x} \\ & =(I+\epsilon J)^{-1}(L+\Delta L_{1})\hat{x}=v.
\end{align*}
Let \(\Delta L=\Delta L_{1}+\Delta L_{2}\).
Then \((L+\Delta L)\hat{x}=v\), with
\begin{align*}
\Delta L=\Delta L_{1}+\Delta L_{2} & =-\epsilon J(I+\epsilon J)^{-1}
L+(I-\epsilon J(I+\epsilon J)^{-1})\Delta L_{1}\\
& =-\epsilon J(I+\epsilon J)^{-1}L+(I+\epsilon J)^{-1}\Delta L_{1}.
\end{align*}
Hence the bound can be presented as
\begin{align*}
\normfro{\Delta L} & \leq\normfro{\epsilon J(I+\epsilon J)^{-1}L}
+\normfro{(I+\epsilon J)^{-1}\Delta L_{1}}
\leq\frac{\beta+\epstrsv(n)}{1-\beta}\normfro{L}.\qedhere
\end{align*}
\end{proof}

\begin{lemma}
\label{lemma:error-unitary}
Let \(H\in\mathbb{C}^{n\times n}\) be a unitary matrix
computed through Householder transformations, \(v\in\mathbb{C}^{n}\).
For applying \(H\) to \(v\) containing a small perturbation, i.e.,
\(u=H\bigl((I+\epsilon J)v\bigr)\), with \(\abs{\epsilon}\leq\beta<1\),
\(\normfro{J}=1\), there exists \(\delta u\) such that
\[
\hat{u}=Hv+\delta u,
\]
where
\[
\norm{\delta u}\leq\frac{\beta+\epsunmqr(n)}{1-\beta}\norm{\hat{u}}.
\]
\end{lemma}

\begin{proof}
From Assumption~\ref{assump:unmqr} we see that
there exists \(\delta u_{1}\) such that
\[
\hat{u}=\fl\bigl(H\bigl((I+\epsilon J)v\bigr)\bigr)
=H(I+\epsilon J)v+\delta u_{1}
\]
where \(\norm{\delta u_{1}}\leq\epsunmqr(n)\cdot\norm{\hat{u}}\).
Let \(\delta u_{2}=\epsilon HJv\) and \(\delta u=\delta u_{1}+\delta u_{2}\).
Then \(\hat{u}=Hv+\delta u\) with
\[
\norm{\delta u}\leq\norm{\delta u_{1}}+\norm{\delta u_{2}}
\leq\epsunmqr(n)\cdot\norm{\hat{u}}+\abs{\epsilon}\cdot\norm{HJv}.
\]
Note that
\[
\norm{HJv}\leq\norm{H}\norm{J}\norm{v}\leq\norm{v}
=\norm{Hv}\leq\norm{\hat{u}}+\norm{\delta u}.
\]
Thus
\[
\norm{\delta u}\leq\epsunmqr(n)\cdot\norm{\hat{u}}
+\beta\bigl(\norm{\hat{u}}+\norm{\delta u}\bigr),
\]
which is equivalent to
\[
\norm{\delta u}\leq\frac{\beta+\epsunmqr(n)}{1-\beta}\norm{\hat{u}}.\qedhere
\]
\end{proof}

\section{The singular values of the preconditioned matrix}
\begin{lemma}
\label{lemma:appendix}
Let
\[
X=\begin{bmatrix}I_{m} & 0 & Z_{1} \\ 0 & 0 & \bar{I}_{p} \\
Z_{1}\herm & \bar{I}_{p}\herm & 0\end{bmatrix},
\]
as in Section~\ref{sec:gmres},
where \(\bar{I}_{p}=[0,I_{p}]\in\mathbb{C}^{p\times n}\). Then
\[
\sigma_{\max}(X)=\lambda_{3},\qquad\sigma_{\min}(X)=\lambda_{2},
\]
where \(\lambda_{1}\approx-1.2470\), \(\lambda_{2}\approx0.4450\),
and \(\lambda_{3}\approx1.8019\) are the roots of
\(\lambda^{3}-\lambda^{2}-2\lambda+1=0\).
\end{lemma}

\begin{proof}
As in Section~\ref{sec:gmres}, we mainly discuss the case
\(m\geq n\geq p\) as the method for the rest is similar.
Let \(\check{I}_{n}=[I_{n},0]\herm\in\mathbb{C}^{m\times n}\).
Note that \(Z_{1}=[Z_{1},Z_{2}][I_{n},0]\herm=Z\check{I}_{n}\).
Therefore
\[
X=\begin{bmatrix}I_{m} & 0 & Z\check{I}_{n} \\ 0 & 0 & \bar{I}_{p} \\
\check{I}_{n}Z\herm & \bar{I}_{p} & 0\end{bmatrix}
=\begin{bmatrix}Z & 0 & 0 \\0 & I_{p} & 0 \\ 0 & 0 & I_{n}\end{bmatrix}
\begin{bmatrix}I_{m} & 0 & \check{I}_{n} \\ 0 & 0 & \bar{I}_{p} \\
\check{I}_{n} & \bar{I}_{p} & 0\end{bmatrix}
\begin{bmatrix}Z\herm & 0 & 0 \\ 0 & I_{p} & 0 \\ 0 & 0 & I_{n}\end{bmatrix}.
\]
Let
\[
Y=\begin{bmatrix}I_{m} & 0 & \check{I}_{n} \\ 0 & 0 & \bar{I}_{p} \\
\check{I}_{n} & \bar{I}_{p} & 0\end{bmatrix}
=\begin{bmatrix}I_{n-p} & 0 & 0 & 0 & I_{n-p} & 0 \\ 0 & I_{p} & 0 & 0 & 0 &
I_{p} \\ 0 & 0 & I_{m-n} & 0 & 0 & 0 \\ 0 & 0 & 0 & 0 & 0 & I_{p} \\
I_{n-p} & 0 & 0 & 0 & 0 & 0 \\ 0 & I_{p} & 0 & I_{p} & 0 & 0\end{bmatrix}.
\]
Then \(X\) is similar to \(Y\) and has the same eigenvalues.

\medskip
A direct calculation gives the characteristic polynomial
\[
\det(\lambda I_{m+n+p}-Y)=(\lambda^{3}-\lambda^{2}-2\lambda+1)^{p}
(\lambda^{2}-\lambda-1)^{n-p}(\lambda-1)^{m-n}.
\]
Therefore
\[
\Lambda(X)=\Lambda(Y)=\begin{cases}\{1,(1\pm\sqrt{5})/2,\lambda_{1},\lambda_{2},
\lambda_{3}\}, & \quad\text{if}\quad m>n>p,\\ \{(1\pm\sqrt{5})/2,\lambda_{1},
\lambda_{2},\lambda_{3}\}, & \quad\text{if}\quad m=n>p,\\
\{1,\lambda_{1},\lambda_{2},\lambda_{3}\}, & \quad\text{if}\quad m>n=p,\\
\{\lambda_{1},\lambda_{2},\lambda_{3}\}, & \quad\text{if}\quad m=n=p.
\end{cases}
\]

Since \(X\) is Hermitian, the singular values of \(X\)
equal to the absolute values of the eigenvalues.
Thus
\[
\sigma_{\max}(X)=\lambda_{3},\qquad\sigma_{\min}(X)=\lambda_{2}.\qedhere
\]
\end{proof}
\medskip
\section*{Declarations}
\begin{itemize}
\item  Competing interests: We declare that the authors have
no competing interests that might be perceived to influence
the results and/or discussion reported in this paper.

\item  Funding: B. Gao and M. Shao are partially supported by
National Key R\&D Program of China under Grant No.\ 2023YFB3001603.
Y. Ma is supported by the European Union (ERC, inEXASCALE, 101075632).
Views and opinions expressed are those of the authors
only and do not necessarily reflect those of
the European Union or the European Research Council.
Neither the European Union nor the granting
authority can be held responsible for them.

\item  Authors' contributions: The authors contributed equally to this work.

\item  Acknowledgments: The authors thank Erin Carson, Ieva
Dau\v{z}ickait\.{e}, and Zhuang-Ao He for helpful discussions.
The authors are also grateful to anonymous reviewers for constructive comments.
\end{itemize}
\end{document}